\documentclass[11pt,reqno,a4paper]{amsart}

\usepackage[T1]{fontenc}
\usepackage{lmodern}
\usepackage[a4paper,left=0.85in,right=0.85in,top=0.9in,bottom=0.9in]{geometry}
\usepackage{amsmath,amsthm,amssymb,mathtools,mathrsfs,enumitem}
\usepackage{xcolor}
\usepackage[expansion=false]{microtype}
\usepackage[hidelinks]{hyperref}
\usepackage{etoolbox}

\renewcommand{\epsilon}{\varepsilon}

\newcommand{\R}{\mathbb R}
\newcommand{\C}{\mathbb C}
\newcommand{\unitdisc}{\mathbb D}
\newcommand{\Aut}{\operatorname{Aut}}
\newcommand{\Ric}{\operatorname{Ric}}
\newcommand{\Lie}{\mathcal L}
\newcommand{\tr}{\operatorname{tr}}
\newcommand{\ddbar}{\sqrt{-1}\,\partial\bar\partial}
\newcommand{\ii}{\sqrt{-1}}
\newcommand{\id}{\operatorname{id}}
\newcommand{\PSH}{\operatorname{PSH}}

\newcommand{\SPSH}{\operatorname{SPSH}}
\newcommand{\dist}{\operatorname{dist}}

\newcommand{\B}{\mathbb B}
\newcommand{\Reg}{\operatorname{Reg}}
\newcommand{\Sing}{\operatorname{Sing}}

\newtheorem{theorem}{Theorem}[section]
\newtheorem{proposition}[theorem]{Proposition}
\newtheorem{lemma}[theorem]{Lemma}
\newtheorem{corollary}[theorem]{Corollary}

\theoremstyle{definition}
\newtheorem{remark}[theorem]{Remark}

\newtheorem{remarks}[theorem]{Remarks}
\AtBeginEnvironment{remarks}{%
  \setlist[enumerate,1]{%
    leftmargin=0pt,
    itemindent=2em,
    labelwidth=1.5em,
    labelsep=0.5em,
    align=left
  }%
}

\title[Expanding solitons and Steinness]%
{Steinness and Complete Expanding K\"ahler--Ricci Solitons}
\author{Soumya Ganguly}
\address{Department of Mathematics, Rutgers University--New Brunswick,
110 Frelinghuysen Rd, Piscataway, NJ 08854, USA}
\email{soumya.ganguly@rutgers.edu}
\date{}
\hypersetup{pdftitle={Steinness and Complete Expanding K\"ahler--Ricci Solitons},
pdfauthor={Soumya Ganguly}}
\begin{document}

\begin{abstract}
We prove that complete expanding K\"ahler--Ricci
solitons with complete soliton fields characterize
Steinness of relatively compact domains in complex
manifolds, provided their closures admit
neighborhoods with holomorphic immersions into
Euclidean space. We construct canonical,
instantaneously complete K\"ahler--Ricci flows
from positive closed $(1,1)$-currents with global
plurisubharmonic potentials and vanishing Lelong
numbers. Nonzero currents homogeneous under a
complete real holomorphic flow yield nontrivial
expanders with prescribed soliton fields. We obtain
a product criterion and uniformization results for
balls and finite-type domains in complex surfaces.
If $D$ is a relatively compact Stein domain in
a reduced pure-dimensional Stein space of dimension
at least two and $\Sing D$ is finite and nonempty,
we prove that $\Aut(\Reg D)$ is compact, excluding
nontrivial expanders with complete soliton fields.
\end{abstract}

\vspace{-5 cm}

\maketitle

\section{Introduction and Main Results}\label{sec:introduction}

Mok--Yau~\cite{MokYau} proved that a bounded
domain in $\C^n$ is pseudoconvex if and only if
it admits a complete K\"ahler--Einstein metric
of negative Einstein constant. For relatively compact
domains in Stein manifolds, the corresponding
characterization of Steinness follows from
\cite[Theorem~2.14]{MokICM2026}.
After rescaling, such a metric $\omega_E$ satisfies
\begin{equation}\label{eq:KE-normalization}
\Ric_{\omega_E}=-K\omega_E,\qquad K>0.
\end{equation}
It gives the self-similar solution
$\omega(t)=(1+Kt)\omega_E$ of the K\"ahler--Ricci flow
\begin{equation}\label{eq:KRF}
\frac{\partial\omega(t)}{\partial t}
=-\Ric_{\omega(t)}.
\end{equation}

More generally, K\"ahler--Ricci solitons are
self-similar solutions evolving by rescaling and
pullback under one-parameter families of
biholomorphisms. In the expanding case, the
soliton equation takes the form
\begin{equation}\label{eq:expanding-soliton}
\Ric_\omega+K\omega=\Lie_Y\omega,
\qquad K>0,
\end{equation}
where $Y$ is a real holomorphic vector field and
$K$ is the \emph{expansion constant}. We require
both the K\"ahler metric $\omega$ and the soliton
field $Y$ to be complete. An expanding soliton is
\emph{trivial} if $\Lie_Y\omega=0$, equivalently if
$\Ric_\omega=-K\omega$, and \emph{nontrivial}
otherwise. Thus a complete real holomorphic
Killing field $Y$ of $\omega_E$, including $Y=0$,
yields a trivial expander. We recall the shrinking,
steady and expanding classification in
Section~\ref{sec:preliminaries}.

We show that Steinness remains necessary when
the K\"ahler--Einstein condition is replaced by
the expanding soliton equation.

\begin{theorem}\label{thm:stein-ambient-expander}
Let $Z$ be a Stein manifold of complex dimension $n\ge1$, and let
$\Omega\Subset Z$ be a domain. Then the following are equivalent:
\begin{enumerate}[label=\textup{(\roman*)}]
\item $\Omega$ is Stein.
\item There exist a complete K\"ahler metric $\omega$ on $\Omega$,
a complete real holomorphic vector field $Y$ on $\Omega$, and a
constant $K>0$ satisfying \eqref{eq:expanding-soliton}.
\end{enumerate}
\end{theorem}

Theorem~\ref{thm:stein-ambient-expander} is a soliton analogue
of Mok--Yau's characterization; the existence direction follows
from their complete negative K\"ahler--Einstein metric.

By the solution of the Levi problem
\cite{Hormander1990}, Theorem~\ref{thm:stein-ambient-expander}
characterizes pseudoconvexity of bounded domains
in $\C^n$. The same characterization extends
to relatively compact domains in arbitrary complex
manifolds whose closures admit neighborhoods
with holomorphic immersions into Euclidean space;
see Remarks~\ref{rem:stein-consequences}.

Ohsawa~\cite{Ohsawa1980} proved that a domain
admitting a complete K\"ahler metric is locally
Stein at its boundary under a $C^1$ boundary
assumption. Diederich--Pflug
\cite[Theorem~1]{DiederichPflug1981}
proved the same conclusion under the condition
$\Omega=\operatorname{int}(\overline{\Omega})$, and consequently Steinness when the ambient
manifold is Stein. Theorem~\ref{thm:stein-ambient-expander}
requires no boundary regularity, but assumes
that the complete metric satisfies the expanding
soliton equation.
Sha~\cite[Theorems~1.2--1.3]{ShaBoundedPseudoconvexKRS2026}
proved K\"ahler--Einstein rigidity for complete
solitons under bounded-geometry and Ricci-pinching
hypotheses when the soliton field has closed metric
dual. He also obtained a ball characterization for
Bergman solitons on smoothly bounded strongly
pseudoconvex domains.

We next construct nontrivial expanders with
prescribed soliton fields. Let $(\omega,Y)$ satisfy \eqref{eq:expanding-soliton}, and let $\{\Phi_s\}_{s\in\R}$ denote the flow generated by $Y$. The associated self-similar K\"ahler--Ricci flow, normalized by
$\omega(1/K)=\omega$, is
\begin{equation}\label{eq:selfsimilar-flow}
\omega(t)
=Kt\,\Phi_{-\frac{1}{K}\log(Kt)}^*\omega,
\qquad t>0.
\end{equation}
This family satisfies
$\partial_t\omega(t)=-\Ric_{\omega(t)}$;
see \cite[Proposition~2.2]{Chowbooksolitonlow}.

Suppose that $\omega(t)$ converges weakly
in the sense of currents to $T$ as
$t\downarrow0$. Since positivity and closedness
are preserved under weak convergence, $T$
is a positive closed $(1,1)$-current.
The self-similarity formula
\eqref{eq:selfsimilar-flow} gives
\[
\Phi_s^*\omega(t)
=e^{Ks}\omega(e^{-Ks}t),
\qquad s\in\R,\quad t>0.
\]
Passing to the weak limit as $t\downarrow0$,
we obtain
\begin{equation}\label{eq:current-homogeneity}
\Phi_s^*T=e^{Ks}T,
\qquad s\in\R.
\end{equation}

The following theorem constructs a complete
nontrivial expander from a nonzero homogeneous
current on a relatively compact Stein domain.

\begin{theorem}\label{thm:current-to-soliton}
Let $\Omega\Subset Z$ be a Stein domain in a
complex manifold. Suppose that a neighborhood of
$\overline\Omega$ admits a holomorphic immersion
into Euclidean space. Let $Y$ be a complete
real holomorphic vector field on $\Omega$
with flow $\{\Phi_s\}_{s\in\R}$. Suppose that,
for some $K>0$, there exists a nonzero
positive closed $(1,1)$-current $T$ on
$\Omega$ satisfying
\eqref{eq:current-homogeneity}.
Then there exists a complete K\"ahler
metric $\omega$ on $\Omega$ satisfying
\eqref{eq:expanding-soliton} and
$\Lie_Y\omega\not\equiv0$.
\end{theorem}

More generally, we establish the existence
of a canonical, instantaneously complete
K\"ahler--Ricci flow from positive closed
$(1,1)$-currents with global plurisubharmonic
potentials and vanishing Lelong numbers on
relatively compact Stein domains whose closures admit
holomorphically immersed neighborhoods
(Theorems~\ref{thm:rough-forward}
and~\ref{thm:instant-complete}).
Our construction permits singular positive closed
$(1,1)$-currents as initial data, going beyond
the smooth semipositive forms considered by
Chau--Lee~\cite[Theorem~1.4]{ChauLee2020}.
See also Ge--Lin--Shen
\cite{KRflowpscvdom2019} for
instantaneously complete flows from smooth
initial metrics. The canonical flow is compatible with biholomorphisms
and constant rescalings
(Proposition~\ref{prop:functoriality}).
In particular, homogeneous initial currents
give rise to self-similar K\"ahler--Ricci flows
(see Proposition~\ref{prop:homogeneous-selfsimilar}).

There is also a natural question about products:
if a relatively compact domain carries a complete
nontrivial expander, when does its product with another
relatively compact domain carry one? Under the
holomorphic immersed-neighborhood hypothesis,
Remarks~\ref{rem:stein-consequences}\textup{(2)}
and the existence of complete negative K\"ahler--Einstein
metrics give the following criterion.

\begin{theorem}\label{thm:product-stability}
Let $D\Subset Z_1$ and $D'\Subset Z_2$ be domains in complex
manifolds, and suppose that neighborhoods of $\overline D$ and
$\overline{D'}$ admit holomorphic immersions into Euclidean spaces.
Suppose $D$ admits a complete nontrivial expanding
K\"ahler--Ricci soliton with complete soliton field $Y_D$ and
expansion constant $K_0>0$. Then $D\times D'$ admits a complete
nontrivial expander with complete soliton field if and only if
$D'$ is Stein. When $D'$ is Stein, one may take the field
$(Y_D,0)$ with expansion constant $K_0$. For every $K>0$ one may
take the field $((K/K_0)Y_D,0)$ with expansion constant $K$.
\end{theorem}

When $D'$ is Stein, the product also admits a
complete K\"ahler--Einstein metric of negative
Einstein constant. Thus $D\times D'$ carries
both trivial and nontrivial expanding solitons.

The noncompactness of a nontrivial soliton flow,
together with the Wong--Rosay theorem for complex
manifolds \cite[Theorem~1]{GaussierKimKrantz2002},
gives the following rigidity result.

\begin{theorem}\label{thm:smooth-ball}
Let $X$ be a complex manifold of dimension $n\ge2$,
and let $\Omega\Subset X$ be a domain whose bounded
holomorphic functions separate points. Suppose that
every boundary orbit accumulation point of
$\Aut(\Omega)$ is a $C^2$ strongly pseudoconvex
boundary point. Then $\Omega$ admits a complete
nontrivial expanding K\"ahler--Ricci soliton with
complete real holomorphic soliton field if and
only if $\Omega$ is biholomorphic to $\B^n$.
\end{theorem}

The point-separation hypothesis holds whenever a
neighborhood of $\overline\Omega$ embeds holomorphically
into some $\C^N$, in particular when $X$ is Stein
\cite[Theorem~3]{NarasimhanEmbedding}.

In complex dimension two, Seshadri--Verma
\cite[Theorem~1.1, p.~1274]{SeshadriVerma2008}
showed that if an automorphism orbit of a domain
in a complex surface accumulates at a boundary
point where the boundary is $C^\infty$ smooth,
pseudoconvex and of finite type in a neighborhood,
then the domain is biholomorphic to a model
of the form
\begin{equation}\label{eq:surface-model}
\mathcal G_H
=\{(z,w)\in\C^2:2\Re w+H(z,\bar z)<0\},
\end{equation}
where $H$ is a nonzero real homogeneous subharmonic
polynomial of even degree without harmonic terms. The maps $\tau_s(z,w)=(z,w+\ii s)$, $s\in\R$,
are automorphisms of $\mathcal G_H$, with infinitesimal
generator $Y_0=\partial/\partial(\Im w)$.
Suppose $\Omega\Subset X$ has a neighborhood of its closure
that embeds holomorphically into $\C^N$. If an automorphism
orbit accumulates at a boundary point where the
boundary is $C^\infty$ smooth, pseudoconvex
and of finite type in a neighborhood,
Corollary~\ref{thm:surface-local} gives
complete nontrivial expanders with soliton fields
$F^*Y_0$ and $-F^*Y_0$ for every biholomorphism
$F:\Omega\to\mathcal G_H$.
If the boundary is smooth, pseudoconvex and of finite
type near every boundary orbit accumulation point,
Corollary~\ref{thm:surface-global} characterizes the
existence of a complete nontrivial expander by
noncompactness of $\Aut(\Omega)$ and biholomorphic
equivalence to some $\mathcal G_H$.

For domains in singular Stein spaces, we prove compactness of the automorphism group of the regular locus without assuming boundary regularity or irreducibility.

\begin{theorem}\label{thm:compact}
Let $Z$ be a reduced Stein space of pure
dimension $n\ge2$, and let $D\Subset Z$
be a connected Stein domain. If $\Sing D$
is finite and nonempty, then
$\Aut(\Reg D)$ is compact.
\end{theorem}

Together with
Lemma~\ref{pcvx:lem:compact-closure-rigidity},
this gives the following consequence.

\begin{corollary}\label{cor:obstruction}
Let $Z$ be a reduced Stein space of pure
dimension $n\ge2$, and let $D\Subset Z$
be a connected Stein domain with $\Sing D$ finite.
If $\Reg D$ admits a complete nontrivial
expanding K\"ahler--Ricci soliton with
complete real holomorphic soliton field,
then $D$ is smooth.
\end{corollary}

\begin{corollary}\label{cor:ball}
Let $Z$ be a reduced Stein space of pure
dimension $n\ge2$, and let $D\Subset Z$
be a connected Stein domain with $\Sing D$
finite. Suppose that every boundary orbit
accumulation point of $\Aut(D)$ belongs to
$\Reg Z$ and is a $C^2$ strongly pseudoconvex
boundary point. Then $\Reg D$ admits a complete
nontrivial expanding K\"ahler--Ricci soliton
with complete real holomorphic soliton field
if and only if $D$ is biholomorphic to $\B^n$.
\end{corollary}

\begin{corollary}\label{thm:bordered-surface}
Let $Z$ be a reduced Stein space of pure
dimension two, and let $D\Subset Z$
be a connected Stein domain satisfying
$\partial_ZD\subset\Reg Z$. Suppose that near every boundary orbit
accumulation point of $\Aut(D)$ the boundary is
$C^\infty$ smooth, pseudoconvex and of finite type.
Then $\Reg D$ admits
a complete nontrivial expanding
K\"ahler--Ricci soliton with complete
real holomorphic soliton field if and
only if $D$ is biholomorphic to a model
$\mathcal G_H$ of the form
\eqref{eq:surface-model}.
In this case, for every biholomorphism
$F:D\to\mathcal G_H$ and every $K>0$,
there exist complete nontrivial expanding
solitons on $D$ with expansion constant $K$
and soliton vector fields $F^*Y_0$ and
$-F^*Y_0$.
\end{corollary}

The uniformization conclusions in
Theorem~\ref{thm:smooth-ball} and
Corollaries~\ref{thm:surface-local},
\ref{thm:surface-global}, \ref{cor:ball}, and~\ref{thm:bordered-surface}
are similar in spirit to the Bergman--Einstein
rigidity and classification results for smoothly
bounded strongly pseudoconvex domains and Stein
spaces with strongly pseudoconvex boundary;
see Huang--Li~\cite{HuangLi2023KEEisBergStein},
Ebenfelt--Xiao--Xu~\cite{EXX2022ChengSteinabelian},
Ganguly--Sinha~\cite{GangulySinha2026ChengStein2d},
and Ganguly--Sahi
\cite{gangulysahi2026classificationsteinspacesbergmaneinstein}.

Section~\ref{sec:preliminaries} fixes our notation and
conventions. Section~\ref{sec:pcvx} proves
Theorem~\ref{thm:stein-ambient-expander}
and Lemma~\ref{pcvx:lem:compact-closure-rigidity}. Section~\ref{sec:analytic-foundations} establishes
the canonical flow and its completeness in
Theorems~\ref{thm:rough-forward}
and~\ref{thm:instant-complete}, and proves
Theorems~\ref{thm:current-to-soliton}
and~\ref{thm:product-stability}.
Section~\ref{sec:boundary-uniformization}
proves Theorem~\ref{thm:smooth-ball}
and Corollaries~\ref{thm:surface-local}
and~\ref{thm:surface-global}.
Section~\ref{sec:bordered-stein-spaces}
proves Theorem~\ref{thm:compact} and its consequences,
Corollaries~\ref{cor:obstruction}, \ref{cor:ball},
and~\ref{thm:bordered-surface}.

\section{Conventions and Preliminaries}
\label{sec:preliminaries}

Throughout the paper, a domain is understood to be a
connected open set. For a subset $A$ of a complex manifold or complex
space $X$, we write $A\Subset X$ if its closure
is compact in $X$. For a domain $\Omega\Subset Z$ we write $\partial_Z\Omega$
for the boundary relative to the ambient manifold or space $Z$.
For a reduced complex space $X$, we denote its regular and
singular loci by $\Reg X$ and $\Sing X$, respectively.

For a K\"ahler manifold $(M,J,g)$, we write
$\omega_g=g(J\cdot,\cdot)$ for its K\"ahler form.
When no confusion is likely, we identify
the metric with its associated form and write
$\omega=\omega_g$. In local holomorphic coordinates,
\[
\omega=\ii g_{j\bar k}\,dz^j\wedge d\bar z^k,
\qquad
\Ric_\omega=-\ddbar\log\det(g_{j\bar k}).
\]
The Riemannian Ricci tensor and the Ricci
form are related by
\[
\Ric_\omega(V,W)=\Ric_g(JV,W)
\]
for all real tangent vectors $V,W$.
We write
\[
\omega_{\mathrm{Euc}}
=\ii\sum_{j=1}^n dz^j\wedge d\bar z^j,
\qquad
dV_{\mathrm E}=\omega_{\mathrm{Euc}}^n.
\]
Thus, for a K\"ahler form $\omega$ on a domain
in $\C^n$, we have
$\omega^n/dV_{\mathrm E}=\det(g_{j\bar k})$.
For a real $(1,1)$-form $\alpha$, we set
\[
\tr_\omega\alpha
=n\frac{\alpha\wedge\omega^{n-1}}{\omega^n}.
\]
Our conventions for the complex Laplacian and
scalar curvature are
\[
\Delta_\omega f
=\tr_\omega(\ddbar f)
=g^{j\bar k}f_{j\bar k},
\qquad
R_\omega=\tr_\omega\Ric_\omega.
\]
If $\Delta_g$ denotes the Riemannian
Laplacian $\tr_g\nabla^2$, and $S_g$
the Riemannian scalar curvature, then
\[
\Delta_g=2\Delta_\omega,
\qquad
S_g=2R_\omega.
\]
For a real-valued smooth function $f$,
the Riemannian gradient and the complex
differential satisfy
\[
|\nabla^g f|_g^2
=2|\partial f|_\omega^2.
\]
The Riemannian volume form is
$dV_g=\omega^n/n!$. We also set
\[
d^c=\ii(\bar\partial-\partial),
\qquad
dd^c=2\ddbar.
\]

For a vector field $Y$ with local flow $\Phi_t$,
we use the convention
\[
\Lie_Y A
=\left.\frac{d}{dt}\right|_{t=0}\Phi_t^*A
\]
for tensor fields $A$, where the pullback of a
vector field $V$ is
$\Phi_t^*V=(\Phi_{-t})_*V$.
In particular, $\Lie_YV=[Y,V]$.
More generally, if $F:M\to N$ is a biholomorphism
and $Y$ is a real vector field on $N$, we write
$F^*Y=(F^{-1})_*Y$.

A real vector field $Y$ on a complex manifold
$(M,J)$ is called \emph{real holomorphic} if
$\Lie_YJ=0$, or equivalently, if its local flow
consists of biholomorphisms. We call $Y$
\emph{complete} if its flow is defined for every
real time. Completeness of a K\"ahler metric
always refers to geodesic completeness.

For a complex manifold or complex space $X$,
we write $\Aut(X)$ for its group of biholomorphic
automorphisms and $\Aut_0(X)$ for the identity
component of $\Aut(X)$. We equip automorphism
and isometry groups with the compact-open
topology, that is,
the topology of uniform convergence on
compact subsets. On the automorphism group of a complex
manifold, Cauchy estimates give smooth convergence on
compact subsets. Continuity of inversion follows by applying
the holomorphic inverse function theorem near each point
of the limiting map.

A holomorphic $\R$-action on a complex
manifold $M$ is a smooth map
$\Phi:\R\times M\to M$ such that
$\Phi_t=\Phi(t,\cdot)$ is holomorphic
for every $t\in\R$ and
\begin{equation}\label{eq:real-action-group-law}
\Phi_0=\operatorname{id}_M,
\qquad
\Phi_{s+t}=\Phi_s\circ\Phi_t
\quad(s,t\in\R).
\end{equation}
In particular, $\Phi_t^{-1}=\Phi_{-t}$,
so each $\Phi_t$ is a biholomorphism.
The infinitesimal generator of the action
is a complete real holomorphic vector field,
and conversely, every complete real
holomorphic vector field generates such
an action. The action is called \emph{free} if
\[
\Phi_t(p)=p
\quad\Longrightarrow\quad t=0
\]
for every $p\in M$. 
It is called
\emph{proper} if the map
\[
\R\times M\longrightarrow M\times M,
\qquad
(t,p)\longmapsto(p,\Phi_t(p))
\]
is proper. More generally, an action of a Lie group
$H$ on $M$ is proper if the map
\[
H\times M\longrightarrow M\times M,
\qquad
(h,p)\longmapsto(p,h\cdot p)
\]
is proper.

For the general soliton equation
$\Ric_\omega=\lambda\omega+\Lie_Y\omega$,
the terms \emph{shrinking}, \emph{steady}, and
\emph{expanding} correspond to $\lambda>0$,
$\lambda=0$, and $\lambda<0$, respectively.
Thus the expansion constant in
\eqref{eq:expanding-soliton} is $K=-\lambda>0$.
For a real holomorphic vector field $Y$,
the equation
$\Ric_\omega=\lambda\omega+\Lie_Y\omega$
is equivalent, in terms of the associated
Riemannian metric $g$, to
\begin{equation}\label{eq:riemannian-soliton}
\Ric_g+\frac12\Lie_{-2Y}g=\lambda g.
\end{equation}

For a complex manifold $M$, we denote by
$\PSH(M)$ the space of plurisubharmonic
functions on $M$. These functions may take the value $-\infty$.
A domain $\Omega\subset\C^n$
is pseudoconvex if it admits a continuous
plurisubharmonic exhaustion, that is,
a continuous function $\rho\in\PSH(\Omega)$
such that
\[
\{\rho<c\}\Subset\Omega
\qquad\text{for every }c\in\R.
\]
By the solution of the Levi problem
\cite[Theorems~2.6.5 and~4.2.8]{Hormander1990}, this is equivalent to Steinness for domains in $\C^n$. On a complex manifold we use the
characterization of Steinness by a smooth strictly
plurisubharmonic exhaustion. We call a smooth real-valued
function $h$ pluriharmonic if $\ddbar h=0$.

For a domain with smooth boundary, pseudoconvexity
and strong pseudoconvexity of the boundary mean
that the Levi form of a local defining function
$\rho$, chosen negative inside the domain, is
respectively positive semidefinite and positive
definite on the complex tangent space.
For smooth boundaries in complex surfaces, finite
type refers to the D'Angelo $1$-type
\cite[Definitions~2.16 and~2.18]{DAngelorealhypordcontc}.
When this type equals $2m$, we say that the
boundary point has type $2m$.

A real $(1,1)$-current $T$ on a complex
manifold of dimension $n$ is called
\emph{positive} if it takes nonnegative
values on smooth compactly supported
strongly positive $(n-1,n-1)$-forms.
It is \emph{closed} if $dT=0$ in the
distributional sense. Every positive closed $(1,1)$-current
has local plurisubharmonic potentials:
near each point there is a
plurisubharmonic function $u$ such that
$T=\ddbar u$. Convergence of $(1,1)$-currents is
understood in the weak sense:
$T_j\rightharpoonup T$ means that
$T_j(\beta)\to T(\beta)$ for every smooth
compactly supported $(n-1,n-1)$-form
$\beta$. For a positive closed $(1,1)$-current $T$
and a point $p$, we denote by $\nu(T,p)$
its Lelong number, with the normalization
of \cite[Chapter~III]{DemaillyCADG}.

For a complex manifold $\Omega$, a complete real
holomorphic vector field $Y$ on $\Omega$ with flow
$\{\Phi_s\}_{s\in\R}$, and a constant $K>0$, we set
\[
\mathcal T_{Y,K}(\Omega)
=
\left\{
T:
\begin{array}{l}
T \text{ is a positive closed real }(1,1)
\text{-current on }\Omega,\\
\Phi_s^*T=e^{Ks}T
\text{ for every }s\in\R
\end{array}
\right\}.
\]
This is a convex cone, which may consist only of
the zero current. On a Stein manifold, a $d$-exact positive closed $(1,1)$-current
has a global plurisubharmonic potential
\cite[Lemma~4.1]{Magnusson2011}. Such a potential is unique up to
a real pluriharmonic function.

Let $\mu$ be a smooth positive volume form on $\Omega$. In holomorphic
coordinates $z=(z^1,\ldots,z^n)$, write $\mu=e^{b_z}dV_z$, where
\[
dV_z=\left(\ii\sum_{j=1}^n dz^j\wedge d\bar z^j\right)^n.
\]

Let $K_\Omega=\Lambda^{n,0}T^*\Omega$ denote
the canonical bundle. The volume form $\mu$
induces a Hermitian metric on $K_\Omega$
for which the local holomorphic frame
$dz^1\wedge\cdots\wedge dz^n$ has squared norm
$e^{-b_z}$. Indeed, under a holomorphic change
of coordinates $z=z(w)$, we have
\[
b_w=b_z+\log\left|\det\frac{\partial z}{\partial w}\right|^2,
\]
which shows that this metric is well-defined.
Its Chern curvature form is $\ddbar b_z$.
We write
\[
\theta_\mu=\ddbar b_z=-\Ric(\mu).
\]
Thus $\theta_\mu$ is a globally defined closed
real $(1,1)$-form.
For every K\"ahler form $\omega$,
\[
-\Ric_\omega=\theta_\mu+\ddbar\log\frac{\omega^n}{\mu}.
\]
Thus the system
\begin{equation}\label{eq:scalar-gauge}
\omega(t)=t\theta_\mu+\ddbar U(t)>0,\qquad
\partial_tU=\log\frac{(t\theta_\mu+\ddbar U)^n}{\mu}
\end{equation}
implies \eqref{eq:KRF}. Moreover, since
$\theta_\mu$ represents $2\pi c_1(K_\Omega)$
and $\ddbar U=\frac12d(d^cU)$ is exact,
\[
[\omega(t)]=t[\theta_\mu]
=2\pi t\,c_1(K_\Omega)
\]
in real de Rham cohomology.

We refer to \eqref{eq:scalar-gauge} as the scalar
potential formulation of the K\"ahler--Ricci flow.
A solution $U(t)$ has initial trace $\psi$ in
$L^1_{\mathrm{loc}}(\Omega)$ if
$U(t)\to\psi$ in $L^1_{\mathrm{loc}}(\Omega)$
as $t\downarrow0$.

If $\widehat\mu=e^b\mu$, where
$b\in C^\infty(\Omega,\R)$, then
$\theta_{\widehat\mu}=\theta_\mu+\ddbar b$.
Consequently, $\widehat U=U-tb$ satisfies
\eqref{eq:scalar-gauge} with $\widehat\mu$
in place of $\mu$, and
\[
t\theta_{\widehat\mu}+\ddbar\widehat U
=t\theta_\mu+\ddbar U.
\]
Since $tb\to0$ locally uniformly as
$t\downarrow0$, this transformation preserves
the initial trace. Similarly, if $p$ is a
time-independent real pluriharmonic function,
then $U+p$ satisfies \eqref{eq:scalar-gauge}
and defines the same K\"ahler forms.
Its initial trace is $\psi+p$ whenever
that of $U$ is $\psi$.

Conversely, suppose that a smooth
K\"ahler--Ricci flow has a global representation
$\omega(t)=t\theta_\mu+\ddbar U(t)$ for $t>0$.
Then
\[
h(t)=\partial_tU(t)-\log\frac{\omega(t)^n}{\mu}
\]
is pluriharmonic in space, since
\eqref{eq:KRF} and the Ricci identity give
$\ddbar h(t)=0$. For any fixed $t_0>0$,
the potential
\[
\widetilde U(t)
=U(t)-\int_{t_0}^t h(s)\,ds
\]
defines the same K\"ahler forms and satisfies
\eqref{eq:scalar-gauge}. The resulting potential need not have the same
initial trace as $U$.

On a domain in $\C^n$, the choice
$\mu=dV_{\mathrm E}$ gives $\theta_\mu=0$.


\section{Complete Expanding Solitons and Steinness}
\label{sec:pcvx}

In this section we prove Theorem~\ref{thm:stein-ambient-expander}.
We begin with a rigidity result for soliton flows with compact
closure. We then establish the metric comparison and analytic-disc
estimates needed to prove that the existence of a complete expanding
soliton implies Steinness.

\subsection{Compact closure of the soliton flow}

Let $(\omega,Y,K)$ satisfy \eqref{eq:expanding-soliton} on a
connected complex manifold $M$, with both $\omega$ and $Y$ complete.
We first record the scalar-curvature estimate used throughout
the converse argument. By
\eqref{eq:riemannian-soliton}, the associated
Riemannian soliton has real dimension $2n$ and
constant $\lambda=-K$. The completeness of $\omega$
and Theorem~1.1(i) of
Catino--Mastrolia--Monticelli--Rigoli
\cite{CatinoMastroliaMonticelliRigoli2016} imply
$S_g\ge-2nK$. Thus
\begin{equation}\label{pcvx:eq:global-scalar-lower}
R_\omega+nK\ge0.
\end{equation}

We next establish rigidity for soliton flows with
relatively compact closure.

\begin{lemma}\label{pcvx:lem:compact-closure-rigidity}
Let $(\omega,Y,K)$ be an expanding soliton on a connected
complex manifold $M$ with both $\omega$ and $Y$ complete,
and let $\{\Phi_s\}_{s\in\R}$ denote the flow of $Y$.
Suppose that the closure of $\{\Phi_s:s\in\R\}$ in $\Aut(M)$
is compact. Then
\[
\Ric_\omega=-K\omega,
\qquad
\Lie_Y\omega=0.
\]
\end{lemma}

\begin{proof}
Taking the trace of
\eqref{eq:expanding-soliton}, we obtain
\[
\begin{aligned}
\Lie_Y(\omega^n)
&=n(\Lie_Y\omega)\wedge\omega^{n-1} =n(\Ric_\omega+K\omega)\wedge\omega^{n-1} =(R_\omega+nK)\omega^n.
\end{aligned}
\]
For each $x\in{M}$, define
$J_x(s)>0$ by
\[
(\Phi_s^*\omega^n)_x
=
J_x(s)\,\omega_x^n.
\]
Differentiating along the flow gives
\[
\frac{d}{ds}\log J_x(s)
=
(R_\omega+nK)(\Phi_s(x))
\ge0.
\]
Thus $J_x(s)$ is nondecreasing
and $J_x(0)=1$. Since the closure of
$\{\Phi_s:s\in\R\}$ in $\Aut({M})$
is compact, we can choose a strictly
increasing sequence of positive integers
$m_j$ such that
$\Phi_{m_j}\to F$ in $\Aut({M})$.
Setting
$t_j=m_{2j}-m_j$,
we have $t_j\to+\infty$ and
$\Phi_{t_j}
=
\Phi_{m_{2j}}\circ\Phi_{m_j}^{-1}
\longrightarrow\id_{M}$.
Since the maps are holomorphic,
their derivatives also converge
uniformly on compact subsets.
Consequently,
\[
J_x(t_j)\longrightarrow1
\qquad(x\in{M}).
\]
For every $s\ge0$, monotonicity
then gives
\[
1=J_x(0)\le J_x(s)
\le\lim_{j\to\infty}J_x(t_j)=1.
\]
Hence $J_x(s)=1$ for all $s\ge0$.
Differentiating at $s=0$ yields
\[
R_\omega(x)+nK=0.
\]
Since $x$ was arbitrary,
we have $R_\omega=-nK$ on ${M}$. Thus $S_g=-2nK$, and
\cite[Theorem~1.1(i)]{CatinoMastroliaMonticelliRigoli2016}
implies that $g$ is Einstein
and $-2Y$ is Killing. Hence $Y$
is Killing, and therefore $\Ric_\omega=-K\omega$, and $\Lie_Y\omega=0$.
\end{proof}
\subsection{Comparison with a complex-hyperbolic metric}
For the following two subsections, let $Z$ be a Stein manifold
and let $\Omega\Subset Z$ be a domain carrying a complete
expanding soliton $(\omega,Y,K)$, with $\omega$ and $Y$ complete.
Let $\{\Phi_s\}_{s\in\R}$ denote the flow of $Y$.
We first establish a lower bound for $\Phi_s^*\omega$
that is uniform in $s$.

\begin{lemma}\label{pcvx:lem:uniform-flow-floor}
There exist a K\"ahler metric $h$ on a neighborhood
of $\overline\Omega$ in $Z$ and a constant $c>0$
such that
\begin{equation}\label{pcvx:eq:uniform-flow-floor}
\Phi_s^*\omega\ge c h
\qquad\text{on }\Omega\quad(s\in\R).
\end{equation}
\end{lemma}

\begin{proof}
Since $Z$ is Stein, the proper holomorphic
embedding theorem~\cite[Theorem~3]{NarasimhanEmbedding}
gives a holomorphic embedding
\[
F:Z\hookrightarrow\C^N.
\]
Choose a Euclidean ball $B\subset\C^N$ with
$F(\overline\Omega)\Subset B$, and let $h_B$
be a complex-hyperbolic metric on $B$
of constant holomorphic sectional curvature
$-\kappa<0$. Set
\[
U:=F^{-1}(B),
\qquad
h:=F^*h_B.
\]
Then $U$ is an open neighborhood of
$\overline\Omega$ in $Z$, and $h$ is a
K\"ahler metric on $U$. By the Gauss equation for K\"ahler submanifolds
\cite[Chapter~IX, Proposition~9.2]{KobayashiNomizuII},
the holomorphic sectional curvature of $h$
does not exceed that of $h_B$. Consequently,
$H_h\le-\kappa$ on $U$.

The metrics $\omega(t)$ defined in
\eqref{eq:selfsimilar-flow} are complete for
every $t>0$, since each is a positive
constant multiple of the pullback of $\omega$
by a biholomorphism of $\Omega$. Fix $t>0$ and $0<\epsilon<t$. The
time-translated family
\[
\widetilde\omega_\epsilon(\tau)
=
\omega(\epsilon+\tau),
\qquad 0\le\tau\le t-\epsilon,
\]
is a complete K\"ahler--Ricci flow on
$[0,t-\epsilon]$. Applying Lee's
parabolic Schwarz estimate
\cite[Lemma~2.1]{L22}, which permits an incomplete
reference metric, with reference metric $h|_\Omega$,
we obtain
\[
\widetilde\omega_\epsilon(\tau)
\ge
\frac{(n+1)\kappa\tau}{2n}\,h,
\qquad 0<\tau\le t-\epsilon.
\]
In particular,
\[
\omega(t)
=
\widetilde\omega_\epsilon(t-\epsilon)
\ge
\frac{(n+1)\kappa(t-\epsilon)}{2n}\,h.
\]
Letting $\epsilon\downarrow0$ yields
\[
\omega(t)
\ge
\frac{(n+1)\kappa t}{2n}\,h,
\qquad t>0.
\]
By \eqref{eq:selfsimilar-flow} and the last estimate,
\[
\Phi_{-K^{-1}\log(Kt)}^*\omega
=\frac{\omega(t)}{Kt}
\ge\frac{(n+1)\kappa}{2nK}\,h,
\qquad t>0.
\]
For any $s\in\R$, taking $t=K^{-1}e^{-Ks}>0$
gives $-K^{-1}\log(Kt)=s$. Hence
\[
\Phi_s^*\omega\ge\frac{(n+1)\kappa}{2nK}\,h,
\qquad s\in\R,
\]
which proves \eqref{pcvx:eq:uniform-flow-floor}
with $c=(n+1)\kappa/(2nK)$.
\end{proof}

\subsection{Analytic discs and their Dirichlet potentials}

We first obtain analytic discs with boundaries in a
fixed compact subset of $\Omega$ and centers
approaching $\partial_Z\Omega$. The discs can be
chosen inside one holomorphic coordinate neighborhood,
with uniformly bounded coordinate derivatives.

\begin{lemma}\label{pcvx:lem:escaping-discs}
If $\Omega\Subset Z$ is not Stein, there exist a
coordinate neighborhood $U$ of a point of
$\partial_Z\Omega$, a holomorphic coordinate map
$z:U\longrightarrow B$ onto a Euclidean ball
$B\subset\C^n$, a point
$p_\infty\in\partial_Z\Omega\cap U$, a compact set
$Q\Subset\Omega\cap U$, and holomorphic maps
$f_j$ from neighborhoods of $\overline\unitdisc$
into $\Omega\cap U$ such that
\begin{equation}\label{pcvx:eq:escaping-discs}
\begin{gathered}
f_j(\overline\unitdisc)\subset\Omega,
\qquad f_j(\partial\unitdisc)\subset Q,\\
f_j(0)\longrightarrow p_\infty,
\qquad
\sup_j\|z\circ f_j\|_{C^1(\overline\unitdisc)}<\infty.
\end{gathered}
\end{equation}
The coordinate map $z$ may be chosen to extend
holomorphically to a neighborhood of $\overline U$.
\end{lemma}

\begin{proof}
Since we assume $\Omega$ is not Stein, by
\cite[\S3, S\"atze~10--11]{DocquierGrauert1960}, there exists
$p\in\partial_Z\Omega$ at which $\Omega$
is not locally Stein. Choose a holomorphic coordinate chart
$z_0:V\to z_0(V)\subset\C^n$ about $p$
and a Euclidean ball $B\Subset z_0(V)$
centered at $z_0(p)$. Set
$U=z_0^{-1}(B)$ and $z=z_0|_U$.
Since $\Omega$ is not locally Stein at $p$,
the set $\Omega\cap U$ is not Stein. At least one connected component
$D$ of $z(\Omega\cap U)$ is then non-Stein,
since a disjoint union of Stein manifolds
is Stein. By the solution of the Levi problem
in $\C^n$ \cite{Hormander1990}, $D$ is not pseudoconvex.

For $\xi\in D$, set
\[
\delta_D(\xi)
=\operatorname{dist}(\xi,\C^n\setminus D),
\qquad
\psi(\xi)=-\log\delta_D(\xi).
\]
Since $D$ is bounded, $\psi$ is a continuous
exhaustion function on $D$. By \cite[Theorem~2.6.7]{Hormander1990},
$\psi$ is not plurisubharmonic.
Consequently, there exist $\xi_0\in D$
and $v\in\C^n$ such that the affine
holomorphic disc
\[
\gamma(\zeta)=\xi_0+\zeta v
\]
satisfies $\gamma(\overline\unitdisc)\Subset D$
and
\[
\psi(\xi_0)>
\frac{1}{2\pi}\int_0^{2\pi}
\psi(\gamma(e^{\ii\theta}))\,d\theta.
\]

Set
\[
a=\psi(\xi_0)-
\frac{1}{2\pi}\int_0^{2\pi}
\psi(\gamma(e^{\ii\theta}))\,d\theta>0.
\]
By the Stone--Weierstrass theorem, choose a real trigonometric polynomial $T_0$ such that
\[
\sup_{\theta\in[0,2\pi]}
|T_0(\theta)-\psi(\gamma(e^{\ii\theta}))|<\frac a4,
\]
and set $T=T_0+a/2$. Then
\[
T(\theta)>\psi(\gamma(e^{\ii\theta}))
\quad(0\le\theta\le2\pi),
\qquad
\frac1{2\pi}\int_0^{2\pi}T(\theta)\,d\theta
<\psi(\xi_0).
\]
Writing
\[
T(\theta)=a_0+\sum_{k=1}^N
\bigl(a_k\cos(k\theta)+b_k\sin(k\theta)\bigr),
\]
define the holomorphic polynomial
\[
P(\zeta)=a_0+\sum_{k=1}^N(a_k-\ii b_k)\zeta^k.
\]
Then $\Re P(e^{\ii\theta})=T(\theta)$ and
\[
\Re P(0)=a_0
=\frac1{2\pi}\int_0^{2\pi}T(\theta)\,d\theta.
\]
Consequently, the mean-value property for $\Re P$
gives
\begin{equation}\label{pcvx:eq:disc-majorant}
\Re P(\zeta)>\psi(\gamma(\zeta))
\quad (|\zeta|=1),
\qquad
\Re P(0)<\psi(\xi_0).
\end{equation}

Since $\delta_D(\xi_0)$ is the distance
from $\xi_0$ to the closed set
$\C^n\setminus D$, there exists
$b\in\C^n$ such that
\[
|b|=\delta_D(\xi_0),
\qquad
\xi_0+b\in\partial D.
\]
For $0\le t\le1$, define
\[
A_t(\zeta)
=\gamma(\zeta)
+t e^{P(0)-P(\zeta)}b.
\]
For $|\zeta|=1$, the inequalities
in \eqref{pcvx:eq:disc-majorant} give
\[
\begin{aligned}
|A_t(\zeta)-\gamma(\zeta)|
&\le\delta_D(\xi_0)
e^{\Re P(0)-\Re P(\zeta)}<e^{-\Re P(\zeta)}
<\delta_D(\gamma(\zeta)).
\end{aligned}
\]
Thus $A_t(\partial\unitdisc)\subset D$
for every $t\in[0,1]$. In particular,
\[
\mathcal K:=
\{A_t(\zeta):0\le t\le1,\ |\zeta|=1\}
\]
is a compact subset of $D$. We have
\[
A_0(\overline\unitdisc)
=\gamma(\overline\unitdisc)\Subset D,
\qquad
A_1(0)=\xi_0+b\in\partial D.
\]
Define
\[
t_*=\inf\{t\in[0,1]:
A_t(\overline\unitdisc)\not\subset D\}.
\]
The compactness of
$A_0(\overline\unitdisc)$ and the continuity
of $(t,\zeta)\mapsto A_t(\zeta)$ imply that
$t_*>0$. Moreover,
$A_t(\overline\unitdisc)\subset D$
for every $t<t_*$. Passing to the limit
as $t\uparrow t_*$ gives
\[
A_{t_*}(\overline\unitdisc)\subset\overline D.
\]
By compactness, the set of parameters $t\in[0,1]$
for which $A_t(\overline\unitdisc)\subset D$
is relatively open. Hence
$A_{t_*}(\overline\unitdisc)\not\subset D$,
and the limiting disc meets $\partial D$. Since
$A_{t_*}(\partial\unitdisc)\subset D$,
there exists $\zeta_*\in\unitdisc$ such that
\[
A_{t_*}(\zeta_*)\in\partial D.
\]

We claim that $A_{t_*}(\zeta_*)\in B$.
Write $B=B(c_B,R)$. Since
$\mathcal K\Subset D\subset B$, there exists
$R_0<R$ such that
$\mathcal K\subset B(c_B,R_0)$.
For each $t\in[0,1]$, the function
$\zeta\mapsto|A_t(\zeta)-c_B|^2$ is
subharmonic. Hence the maximum principle
gives
\[
A_t(\overline\unitdisc)
\subset\overline{B(c_B,R_0)}\Subset B,
\qquad 0\le t\le1.
\]
In particular,
$A_{t_*}(\zeta_*)\in\partial D\cap B$.
Since $D$ is a connected component of
$z(\Omega\cap U)$, every point of
$\partial D\cap B$ belongs to
$z(\partial_Z\Omega\cap U)$.
It follows that
\[
p_\infty:=
z^{-1}(A_{t_*}(\zeta_*))
\in\partial_Z\Omega\cap U.
\]

Choose an automorphism $\varphi$ of $\unitdisc$ satisfying
$\varphi(0)=\zeta_*$. Recall that
$\mathcal K=\bigcup_{0\le t\le1}A_t(\partial\unitdisc)\Subset D$.
Set $\delta=\operatorname{dist}(\mathcal K,\C^n\setminus D)>0$,
and let $Q_0$ be the closed $\delta/2$-neighborhood of
$\mathcal K$, so that $Q_0\Subset D$.

Since $\varphi$ extends continuously to $\overline\unitdisc$,
the map $(t,r,\zeta)\mapsto A_t(\varphi(r\zeta))$ is uniformly
continuous on $[0,1]\times[1/2,1]\times\partial\unitdisc$.
As $\varphi(\partial\unitdisc)=\partial\unitdisc$, we may
choose $\rho\in(1/2,1)$ such that
\[
\sup_{\substack{0\le t\le1\\|\zeta|=1}}
|A_t(\varphi(\rho\zeta))-A_t(\varphi(\zeta))|
<\frac{\delta}{2}.
\]
Consequently,
$A_t(\varphi(\rho\partial\unitdisc))\subset Q_0$
for every $t\in[0,1]$. Choose $t_j\uparrow t_*$ with $0<t_j<t_*$, and define
\[
f_j(\zeta)=
z^{-1}\bigl(A_{t_j}(\varphi(\rho\zeta))\bigr).
\]
Since $t_j<t_*$, the image of
$\overline\unitdisc$ under $f_j$ lies
in $\Omega\cap U$. Furthermore,
\[
f_j(\partial\unitdisc)\subset Q,
\qquad Q:=z^{-1}(Q_0)\Subset\Omega\cap U,
\]
and
\[
f_j(0)
=z^{-1}(A_{t_j}(\zeta_*))
\longrightarrow p_\infty.
\]

Choose $\eta>0$ such that $\rho(1+\eta)<1$.
Since $A_{t_j}(\overline\unitdisc)\subset D$, the maps
\[
z\circ f_j=A_{t_j}\circ\varphi(\rho\,\cdot)
\]
are holomorphic on $(1+\eta)\unitdisc$ and take values
in the bounded set $D\subset B$. Applying Cauchy's
estimates to their coordinate functions on discs of
radius $\eta/2$ centered at points of $\overline\unitdisc$
gives
\[
\sup_j\|z\circ f_j\|_{C^1(\overline\unitdisc)}<\infty.
\] 
This completes the proof.
\end{proof}

For a holomorphic map $f$ defined on a neighborhood
of $\overline\unitdisc$, with
$f(\overline\unitdisc)\subset\Omega$,
let $u_f\in C^\infty(\overline\unitdisc)$
be the unique solution of
\begin{equation}\label{pcvx:eq:zero-dirichlet}
\ddbar u_f=f^*\omega
\quad\text{on }\unitdisc,
\qquad
u_f|_{\partial\unitdisc}=0.
\end{equation}
Since $f^*\omega\ge0$, the maximum principle
gives $u_f\le0$. We record a particular form of
\cite[Theorem~1]{BrezisMerle1991}.

\begin{lemma}\label{pcvx:lem:green-exponential}
Let $f_j$ be holomorphic maps defined on
neighborhoods of $\overline\unitdisc$, with
$f_j(\overline\unitdisc)\subset\Omega$, and set
$u_j=u_{f_j}$. Suppose that
$\int_{\unitdisc}f_j^*\omega\le M$
for some $M>0$ independent of $j$.
Then, for every $0<q<2\pi/M$,
\begin{equation}\label{pcvx:eq:green-exponential}
\sup_j\int_{\unitdisc}e^{-qu_j}\,dA<\infty,
\end{equation}
where $dA$ denotes Euclidean area measure.
\end{lemma}

\begin{proof}
Let $d\mu_j=f_j^*\omega$ and
$0 \leq m_j=\mu_j(\unitdisc)\le M$.
For $a\in\unitdisc$, set
\[
G(\zeta,a)
=\log\left|\frac{\zeta-a}{1-\bar a\zeta}\right|.
\]
For $|\zeta|=1$, we have
$|\zeta-a|=|1-\bar a\zeta|$, so
$G(\,\cdot\,,a)$ vanishes on $\partial\unitdisc$.
Moreover, $\log|1-\bar a\zeta|$ is harmonic
on $\unitdisc$, while
$\Delta_\zeta\log|\zeta-a|=2\pi\delta_a$
in the sense of distributions. Here $\Delta$
is the Euclidean Laplacian and $\delta_a$
is the Dirac mass at $a$.
Since $\ddbar u=\frac12(\Delta u)\,dA$,
we have $\ddbar_\zeta G(\zeta,a)=\pi\delta_a$.
By the Green-function representation for the
Dirichlet problem~\cite[\S2.2.4, Theorem~12]{PDEbookEvans_1998},
\[
u_j(\zeta)=\frac1\pi
\int_{\unitdisc}G(\zeta,a)\,d\mu_j(a).
\]

If $m_j=0$, then $u_j=0$. Otherwise,
$d\nu_j=m_j^{-1}d\mu_j$ is a probability measure,
and Jensen's inequality gives
\[
e^{-qu_j(\zeta)}
\le\int_{\unitdisc}
\left|\frac{1-\bar a\zeta}{\zeta-a}\right|^{qm_j/\pi}
\,d\nu_j(a).
\]
Set $\beta=qM/\pi<2$. Since
\[
|1-\bar a\zeta|^2-|\zeta-a|^2
=(1-|a|^2)(1-|\zeta|^2)\ge0,
\]
the ratio in the preceding integral is at least
$1$. Using $qm_j/\pi\le\beta$ and
$|1-\bar a\zeta|\le2$, Tonelli's theorem gives
\[
\begin{aligned}
\int_{\unitdisc}e^{-qu_j}\,dA
&\le 2^\beta\sup_{a\in\unitdisc}
\int_{\unitdisc}|\zeta-a|^{-\beta}\,dA(\zeta)\\
&\le 2^\beta\int_{|w|<2}|w|^{-\beta}\,dA(w)\\
&=\frac{8\pi}{2-\beta}<\infty.
\end{aligned}
\]
The bound is independent of $j$.
\end{proof}

We now deduce the following consequence of the
Hopf--Rinow theorem, which will be applied to
the sequence $f_j(0)\to p_\infty\in\partial_Z\Omega$
constructed in Lemma~\ref{pcvx:lem:escaping-discs}.

\begin{lemma}\label{lem:complete-ambient-escape}
Let $(M,g)$ be a connected complete Riemannian manifold embedded in a Hausdorff manifold $E$. Suppose that a sequence $p_j\in M$ converges in $E$ to a point $p_\infty\in E\setminus M$. Then, for every nonempty compact set $L\subset M$,
\[
d_g(p_j,L)\longrightarrow+\infty.
\]
Consequently, the lengths of any piecewise smooth curves in $M$ joining $p_j$ to points of $L$ tend to infinity.
\end{lemma}

\begin{proof}
Suppose that $d_g(p_j,L)$ does not tend to infinity. After passing to a subsequence, there exists $R>0$ such that $d_g(p_j,L)\le R$ for every $j$. By the Hopf--Rinow theorem, the closed $R$-neighborhood of the compact set $L$ is compact in $M$. Thus, after passing to a further subsequence, $p_j$ converges to some $p\in M$. Since the inclusion $M\hookrightarrow E$ is continuous, this subsequence also converges to $p$ in $E$, contradicting the uniqueness of limits in the Hausdorff manifold $E$, because $p_j\to p_\infty\notin M$.

The final assertion follows since every piecewise smooth curve joining $p_j$ to $L$ has length at least $d_g(p_j,L)$.
\end{proof}

\subsection{Proof of the Steinness characterization}

\begin{proof}[Proof of Theorem~\ref{thm:stein-ambient-expander}]
We first prove \textup{(i)}$\Rightarrow$\textup{(ii)}. Suppose that $\Omega$ is Stein. If $n\ge2$, it is a domain
of holomorphy in $Z$ by
Docquier--Grauert~\cite[\S3, p.~116]{DocquierGrauert1960}.
By \cite[Theorem~2.14]{MokICM2026}, it admits a complete
K\"ahler--Einstein metric with negative Einstein constant;
see also \cite[\S3.2]{MokYau} for the bounded-domain case.

When $n=1$, let $F:Z\hookrightarrow\C^N$ be a proper
holomorphic embedding
\cite[Theorem~3]{NarasimhanEmbedding}.
Since $F$ is an immersion, at least one of its coordinate
functions is nonconstant on $\Omega$. This function is
bounded on $\Omega$, since $\overline\Omega$ is compact in $Z$.
It lifts to a bounded nonconstant holomorphic function
on the universal cover of $\Omega$.
By the uniformization theorem
\cite[Theorem~8.23]{Schlag2014}, this cover is
biholomorphic to $\mathbb P^1$, $\C$, or $\unitdisc$.
The first two possibilities are excluded by the maximum
principle and Liouville's theorem, respectively.
Thus the universal cover is biholomorphic to $\unitdisc$.
The deck transformations preserve the Poincar\'e metric,
which therefore descends to a K\"ahler--Einstein metric
on $\Omega$ with negative Einstein constant.
The descended metric is complete, since its geodesics
lift to geodesics of the complete Poincar\'e metric.

In either case, write $\Ric_{\omega_0}=-c_{\mathrm{KE}}\omega_0$,
where $\omega_0$ is complete and $c_{\mathrm{KE}}>0$.
For any $K>0$, the rescaled form
$(c_{\mathrm{KE}}/K)\omega_0$
satisfies \eqref{eq:KE-normalization}, since its Ricci form
is unchanged by constant rescaling.
Taking $Y=0$ proves \textup{(ii)}.

We now prove \textup{(ii)}$\Rightarrow$\textup{(i)}.
Suppose that $\Omega$ admits a complete K\"ahler metric $\omega$,
a complete real holomorphic vector field $Y$, and a constant
$K>0$ satisfying \eqref{eq:expanding-soliton}.
Let $\{\Phi_s\}_{s\in\R}$ denote the flow of $Y$.
When $n=1$, $\Omega$ is a connected noncompact
Riemann surface and hence is Stein by
Behnke--Stein~\cite{BehnkeStein1949}.
We may therefore assume $n\ge2$.
Suppose that $\Omega$ is not Stein.
Let $U$, $z$, $f_j$, $Q$, and $p_\infty$ be as in
Lemma~\ref{pcvx:lem:escaping-discs}.
On $\Omega\cap U$, write in the holomorphic
coordinates $z=(z^1,\ldots,z^n)$
\[
\omega_s=\Phi_s^*\omega
=\ii(g_s)_{a\bar b}\,dz^a\wedge d\bar z^b,
\qquad
r_s=\log\det((g_s)_{a\bar b}),
\]
and set $r=r_0$. Since $\Phi_s$ is biholomorphic,
$\Ric_{\omega_s}=\Phi_s^*\Ric_\omega=-\ddbar r_s$.
Moreover, $\partial_s\omega_s=\Phi_s^*(\Lie_Y\omega)$.
Using \eqref{eq:expanding-soliton} and
the identity
$\partial_s\log\det(g_s)
=\tr(g_s^{-1}\partial_s g_s)$,
we obtain
\begin{equation}\label{pcvx:eq:pullback-soliton}
\begin{aligned}
\partial_s\omega_s
&=K\omega_s-\ddbar r_s,\\
\partial_s r_s
&=\tr_{\omega_s}(\partial_s\omega_s)
=(R_\omega+nK)\circ\Phi_s\ge0,
\end{aligned}
\end{equation}
where the last inequality follows from
\eqref{pcvx:eq:global-scalar-lower}.

We first bound the $\omega$-areas of the discs.
As $\Ric_\omega=-\ddbar r$ and
$\Ric_\omega+K\omega=\Lie_Y\omega$,
Cartan's formula on the closed form $\omega$
gives $\Lie_Y\omega=d(\iota_Y\omega)$.
By the convention $dd^c=2\ddbar$,
the soliton equation becomes
\begin{equation}\label{pcvx:eq:DG-local-exactness}
K\omega
=d\left(\iota_Y\omega+\frac12d^c r\right).
\end{equation}
Set
\[
\alpha=\iota_Y\omega+\frac12d^c r.
\]
The one-form $\alpha$ is smooth on $\Omega\cap U$.
Since each $f_j(\overline\unitdisc)$ lies in
$\Omega\cap U$, by Stokes' theorem 
\begin{equation}\label{pcvx:eq:DG-disc-area-Stokes}
K\int_{\unitdisc}f_j^*\omega
=\int_{\partial\unitdisc}f_j^*\alpha.
\end{equation}
The images $f_j(\partial\unitdisc)$ lie in
$Q\Subset\Omega\cap U$, where $\alpha$ is bounded,
and the coordinate derivatives are uniformly bounded
by \eqref{pcvx:eq:escaping-discs}. Thus there is
$M>0$, independent of $j$, such that
\begin{equation}\label{pcvx:eq:disc-area-uniform}
\int_{\unitdisc}f_j^*\omega\le M.
\end{equation}

For each $j$ and $s\in\R$, let
$u_{j,s}\in C^\infty(\overline\unitdisc)$
be the unique solution of
\[
\ddbar u_{j,s}=f_j^*\omega_s
\quad\text{on }\unitdisc,
\qquad
u_{j,s}|_{\partial\unitdisc}=0.
\]
Since $f_j^*\omega_s\ge0$, the maximum principle
gives $u_{j,s}\le0$. Moreover,
$u_{j,0}=u_{f_j}$ by uniqueness of the
Dirichlet problem \eqref{pcvx:eq:zero-dirichlet}.

For each fixed $j$, write
$f_j^*\omega_s=F_j(s,a)\,dA(a)$.
Since $F_j$ is smooth on
$\R\times\overline\unitdisc$ and
\[
\sup_{\zeta\in\overline\unitdisc}
\int_{\unitdisc}|G(\zeta,a)|\,dA(a)
\le\int_{|w|<2}\log(2/|w|)\,dA(w)<\infty,
\]
we may differentiate under the integral sign in
\[
u_{j,s}(\zeta)=\frac1\pi
\int_{\unitdisc}G(\zeta,a)F_j(s,a)\,dA(a)
\]
uniformly for $\zeta\in\overline\unitdisc$
and $s$ in compact intervals. 

By \eqref{pcvx:eq:pullback-soliton},
\[
\begin{aligned}
\ddbar(\partial_su_{j,s})
&=f_j^*(\partial_s\omega_s) =K\ddbar u_{j,s}-\ddbar(r_s\circ f_j).
\end{aligned}
\]
Thus
$\partial_su_{j,s}-Ku_{j,s}+r_s\circ f_j$
is harmonic on $\unitdisc$. Since
$u_{j,s}|_{\partial\unitdisc}=0$ for every $s$,
we also have
$(\partial_su_{j,s})|_{\partial\unitdisc}=0$.
Let $H_{j,s}$ denote the harmonic extension
of $(r_s\circ f_j)|_{\partial\unitdisc}$.
Uniqueness of the harmonic Dirichlet problem
then gives
\begin{equation}\label{pcvx:eq:dirichlet-evolution}
\partial_su_{j,s}
=Ku_{j,s}-r_s\circ f_j+H_{j,s}.
\end{equation}

Fix $S>0$. Since $r_s$ depends continuously
on $(s,p)\in[-S,0]\times Q$, the constant
\[
C_S:=\max_{\substack{-S\le s\le0\\p\in Q}}r_s(p)
\]
is finite. As $f_j(\partial\unitdisc)\subset Q$,
the boundary values of $H_{j,s}$ are bounded above
by $C_S$ for every $j$ and $-S\le s\le0$.
The maximum principle therefore gives
\[
H_{j,s}\le C_S
\quad\text{on }\unitdisc,
\qquad j\ge1,\ -S\le s\le0.
\]
Multiplying \eqref{pcvx:eq:dirichlet-evolution}
by $e^{-Ks}$ and integrating from $-S$ to $0$
gives
\[
\begin{aligned}
u_{j,0}
&=e^{KS}u_{j,-S}
-\int_{-S}^0e^{-Ks}r_s\circ f_j\,ds +\int_{-S}^0e^{-Ks}H_{j,s}\,ds.
\end{aligned}
\]
We know $u_{j,-S}\le0$ and from \eqref{pcvx:eq:pullback-soliton},
$r_s\ge r_{-S}$ for $-S\le s\le0$. Hence 
\begin{equation}\label{pcvx:eq:determinant-upper}
r_{-S}(f_j(\zeta))
\le C_S-\frac{u_{j,0}(\zeta)}{b_S},
\qquad
b_S=\int_{-S}^0e^{-Ks}\,ds
=\frac{e^{KS}-1}{K}.
\end{equation}

By Lemma~\ref{pcvx:lem:uniform-flow-floor},
$\omega_{-S}\ge c h$ on $\Omega$, with $c$
independent of $S$. Since $z$ extends as a coordinate chart
across $\overline U$ and $h$ is a K\"ahler
metric near $\overline\Omega$, compactness
of $\overline\Omega\cap\overline U$ gives
$c_U>0$ such that
\[
h\ge c_U\omega_{\mathrm{Euc}}
\qquad\text{on }\Omega\cap U,
\]
where $\omega_{\mathrm{Euc}}=
\ii\sum_a dz^a\wedge d\bar z^a$ in the
coordinates $z$. Put $c_0=cc_U>0$. Let $0<\lambda_1\le\cdots\le\lambda_n$
denote the eigenvalues of $g_{-S}$
relative to the Euclidean metric. Since $\omega_{-S}\ge c_0\omega_{\mathrm{Euc}}$,
we have $\lambda_a\ge c_0$ for every
$a=1,\ldots,n$. Consequently,
\[
\lambda_n \le c_0^{-(n-1)}\det(g_{-S})
=c_0^{-(n-1)}e^{r_{-S}}.
\]
By \eqref{pcvx:eq:escaping-discs}, there exists
$L>0$, independent of $j$, such that
$|d(z\circ f_j)|_{\mathrm{Euc}}\le L$ on
$\overline\unitdisc$. The preceding eigenvalue
estimate and \eqref{pcvx:eq:determinant-upper}
therefore give
\[
|df_j(\zeta)|_{\omega_{-S}}^2
\le L^2c_0^{-(n-1)}e^{C_S}
\exp\left(-\frac{u_{j,0}(\zeta)}{b_S}\right).
\]
Consequently,
\begin{equation}\label{pcvx:eq:metric-derivative-upper}
|df_j(\zeta)|_{\omega_{-S}}
\le A_S
\exp\left(-\frac{u_{j,0}(\zeta)}{2b_S}\right),
\end{equation}
where $A_S>0$ depends on $S$ but not on $j$.

Fix $\ell>2$. Since $b_S\to\infty$ as
$S\to\infty$, we may choose $S>0$ sufficiently
large that
\[
q:=\frac{\ell}{2b_S}<\frac{2\pi}{M}.
\]
Fix such an $S$ and set
\[
\rho_j(\zeta)
=\exp\left(-\frac{u_{j,0}(\zeta)}{2b_S}\right).
\]
By Lemma~\ref{pcvx:lem:green-exponential}
and \eqref{pcvx:eq:disc-area-uniform},
\[
\sup_j\|\rho_j\|_{L^\ell(\unitdisc)}<\infty.
\]
Let $\ell'=\ell/(\ell-1)\in(1,2)$. With $dA$ denoting
Euclidean area measure, H\"older's inequality gives
\[
\begin{aligned}
\int_0^{2\pi}\!\int_0^1
\rho_j(te^{\ii\theta})\,dt\,d\theta
&=\int_{\unitdisc}
\frac{\rho_j(\zeta)}{|\zeta|}\,dA(\zeta) \le\|\rho_j\|_{L^\ell(\unitdisc)}
\left(\int_{\unitdisc}|\zeta|^{-\ell'}\,dA\right)^{1/\ell'} \le C < \infty,
\end{aligned}
\]
where $C$ is independent of $j$, since $\ell'<2$.
Consequently, for every $j$ there is
$\theta_j\in[0,2\pi)$ such that
\[
\int_0^1\rho_j(te^{\ii\theta_j})\,dt
\le\frac{C}{2\pi}.
\]
For $\gamma_j(t)=te^{\ii\theta_j}$,
$0\le t\le1$, inequality
\eqref{pcvx:eq:metric-derivative-upper} yields
\begin{equation}\label{pcvx:eq:radial-length}
L_{\omega_{-S}}(f_j\circ\gamma_j)
\le\frac{A_SC}{2\pi}
\qquad(j\ge1).
\end{equation}

The metric $\omega_{-S}=\Phi_{-S}^*\omega$
is complete because $\Phi_{-S}$ is a
biholomorphism of $\Omega$ and $\omega$ is
complete. But $f_j(0)\to p_\infty\in\partial_Z\Omega$
and $f_j(e^{\ii\theta_j})\in Q$.
Applying Lemma~\ref{lem:complete-ambient-escape}
to the inclusion $\Omega\hookrightarrow Z$
with the complete metric $\omega_{-S}$ gives
\[
d_{\omega_{-S}}(f_j(0),Q)\longrightarrow+\infty.
\]
This contradicts the uniform length bound
\eqref{pcvx:eq:radial-length}. Therefore
$\Omega$ is Stein.
\end{proof}

\begin{remarks}\label{rem:stein-consequences}
\begin{enumerate}[label=\textup{(\arabic*)}]
\item For bounded domains in $\C^n$,
Theorem~\ref{thm:stein-ambient-expander}
characterizes pseudoconvexity, by the solution
of the Levi problem \cite{Hormander1990}.

\item The equivalence in
Theorem~\ref{thm:stein-ambient-expander}
continues to hold for a domain
$\Omega\Subset Z$ in an arbitrary complex
manifold $Z$, provided a neighborhood of
$\overline\Omega$ admits a holomorphic
immersion into $\C^N$.

To see this, choose
$\overline\Omega\subset X\Subset V\subset Z$
and a holomorphic immersion
$\iota:V\to\C^N$. Let $B\subset\C^N$ be
a ball containing $\iota(\overline X)$,
and let $h_B$ be its complex-hyperbolic
metric. Set $\varphi=|\iota|^2$,
$\beta=\ddbar\varphi$, and $h=\iota^*h_B$.
The immersion hypothesis implies that
$\beta$ and $h$ are K\"ahler metrics on
$X$. By the K\"ahler Gauss equation,
$H_h\le-\kappa<0$, and $h\ge c\beta$
for some positive constants $\kappa,c$.

If $\Omega$ is Stein, it is locally Stein
over $X$: for every Stein coordinate
neighborhood $U\subset X$, the diagonal
of $\Omega\cap U$ is a closed complex
submanifold of the Stein manifold
$\Omega\times U$. Thus Mok--Yau
\cite[\S3.1, case~(iii)]{MokYau}
gives a complete negative K\"ahler--Einstein
metric on $\Omega$, which yields the
required expander with $Y=0$.

Conversely, suppose that $\Omega$ admits
a complete expanding soliton with complete
soliton field. In dimension one, $\Omega$
is noncompact by the immersion hypothesis
and hence Stein by Behnke--Stein
\cite{BehnkeStein1949}. For $n\ge2$,
the proof of
Lemma~\ref{pcvx:lem:uniform-flow-floor},
with $h$ as reference metric, gives
$\Phi_s^*\omega\ge ch$ for all $s\in\R$.
If $\Omega$ were not locally Stein at a
boundary point, the local disc construction
in Lemma~\ref{pcvx:lem:escaping-discs}
and the estimates in the proof of
Theorem~\ref{thm:stein-ambient-expander}
would contradict completeness. Hence
$\Omega$ is locally Stein over $X$.
Since $\varphi$ is strictly plurisubharmonic
on $X$, Elencwajg's theorem
\cite[Th\'eor\`eme~II]{Elencwajg1975}
implies that $\Omega$ is Stein.

\item The ambient hypothesis cannot be omitted, even
for a K\"ahler manifold. Let $S$ be a compact Riemann
surface of genus at least two, $X=S\times\C$, and
$\Omega=S\times\unitdisc\Subset X$. The product of
a K\"ahler metric on $S$ and the Euclidean metric on
$\C$ makes $X$ K\"ahler. By uniformization
\cite[Theorem~8.23]{Schlag2014}, $S$ and $\unitdisc$
carry complete hyperbolic metrics, which may be
scaled to have the same negative Einstein constant
$-K$. Their product is a complete K\"ahler--Einstein
metric on $\Omega$, hence gives a complete expander
with $Y=0$. Nevertheless, $\Omega$ is not Stein:
every holomorphic function on $\Omega$ is constant
on the compact complex curve $S\times\{0\}$,
contradicting holomorphic separability
\cite[p.~328]{Narasimhan1962}.
\end{enumerate}
\end{remarks}

\section{From Positive Closed Currents to Complete Expanding Solitons}
\label{sec:analytic-foundations}

Throughout this section, let $Z$ be a complex manifold,
let $\Omega\Subset Z$ be a Stein domain, and suppose that
a neighborhood of $\overline\Omega$ admits a holomorphic
immersion into Euclidean space. Fix $K>0$.
We first construct K\"ahler--Ricci flows starting from
positive closed $(1,1)$-currents with global plurisubharmonic
potentials and vanishing Lelong numbers. We prove that these flows are instantaneously complete.
For currents satisfying \eqref{eq:current-homogeneity},
we establish self-similarity. When the current
is nonzero, this construction gives a complete
nontrivial expanding K\"ahler--Ricci soliton.

By Remarks~\ref{rem:stein-consequences}\textup{(2)},
$\Omega$ admits a complete K\"ahler--Einstein metric
$\omega_E$ satisfying \eqref{eq:KE-normalization}.

\subsection{The scalar K\"ahler--Ricci flow and comparison}
\label{sec:scalar-comparison}

Choose open neighborhoods
$\overline\Omega\subset X\Subset V\subset Z$ and a
holomorphic immersion $\iota:V\to\C^N$.
Let $B\subset\C^N$ be a Euclidean ball containing
$\iota(\overline X)$ compactly, and let $h_B$
be its complex-hyperbolic metric. Set
\[
h=\iota^*h_B\qquad\text{on }X.
\]
By the Gauss equation for K\"ahler submanifolds
\cite[Chapter~IX, \S9, Propositions~9.1--9.2]{KobayashiNomizuII},
the holomorphic bisectional curvature of $h$ is bounded
above by that of $h_B$, and hence by a negative constant.
In particular, $\Ric_h<0$, and there exists $\kappa>0$
such that the holomorphic sectional curvature of $h$
is at most $-\kappa$.

Set
\[
\mu=h^n,\qquad \theta=-\Ric_h>0,\qquad
\varphi=|\iota|^2,\qquad \beta=\ddbar\varphi.
\]
Both $h$ and $\beta$ are K\"ahler forms on $X$.
Since $\iota(\overline X)\Subset B$, there exists $c>0$
such that $h\ge c\beta$ on $X$.

We will use the following comparison principle
for the scalar equation \eqref{eq:scalar-gauge}.

\begin{lemma}\label{lem:scalar-comparison}

Let $D\Subset Z$ be a domain, let $\mu$ be a smooth
positive volume form on $D$, and let $\tau>0$.
Suppose that $U,V\in C^\infty(D\times(0,\tau])$
and that
\[
A_U=t\theta_\mu+\ddbar U>0,
\qquad
A_V=t\theta_\mu+\ddbar V>0.
\]
Assume that
\begin{align}
\partial_t U&\ge\log\frac{A_U^n}{\mu},\label{eq:comparison-super}\\
\partial_t V&\le\log\frac{A_V^n}{\mu},\label{eq:comparison-sub}
\end{align}
and
\begin{equation}\label{eq:comparison-initial}
\liminf_{t\downarrow0}\inf_D(U(t)-V(t))\ge0.
\end{equation}
If, for every compact interval $I\subset(0,\tau]$,
\begin{equation}\label{eq:comparison-boundary}
\liminf_{z\to\partial D}\inf_{t\in I}(U(t,z)-V(t,z))\ge0,
\end{equation}
then $U\ge V$ on $D\times(0,\tau]$.
\end{lemma}

\begin{proof}
Set $w=U-V$. In local holomorphic coordinates, let
$G_U$ and $G_V$ denote the coefficient matrices of
$A_U$ and $A_V$, respectively. For $0\le\lambda\le1$,
define
\[
G_\lambda=(1-\lambda)G_V+\lambda G_U.
\]
Since $A_U-A_V=\ddbar w$, we have
$\partial_\lambda G_\lambda=(w_{i\bar j})$.
The matrix differentiation formula
\[
\frac{d}{d\lambda}\log\det G_\lambda
=\tr\bigl(G_\lambda^{-1}(w_{i\bar j})\bigr)
\]
and the fundamental theorem of calculus give
\[
\log\frac{A_U^n}{A_V^n}
=a^{i\bar j}w_{i\bar j},
\qquad
(a^{i\bar j})
=\int_0^1G_\lambda^{-1}\,d\lambda.
\]
Since $A_U$ and $A_V$ are positive definite,
so is $(a^{i\bar j})$. Consequently,
\[
\partial_tw\ge a^{i\bar j}w_{i\bar j}.
\]

Fix $t_0\in(0,\tau]$ and $\delta>0$.
By \eqref{eq:comparison-initial}, we may
choose $t_1\in(0,t_0)$ such that
\[
w(t_1,z)>-\delta
\qquad\text{for every }z\in D.
\]
By \eqref{eq:comparison-boundary}, applied
with $I=[t_1,t_0]$, there exists a compact
set $L\Subset D$ such that
\[
w(t,z)>-\delta
\qquad
(z\in D\setminus L,\ t_1\le t\le t_0).
\]
Choose a smoothly bounded relatively compact
domain $D'\subset D$ containing $L$. Then
\[
w(t,z)>-\delta
\qquad
(z\in D\setminus D',\ t_1\le t\le t_0).
\]

Set $v(t,z)=w(t,z)+\delta+\delta t$.
By the choices of $t_1$ and $D'$, we have
\[
v(t_1,z)>0
\quad(z\in\overline{D'}),
\qquad
v(t,z)>0
\quad(z\in\partial D',\ t\in[t_1,t_0]).
\]
Moreover,
\[
\partial_tv\ge a^{i\bar j}v_{i\bar j}+\delta.
\]
Suppose that $v$ were nonpositive somewhere
on $\overline{D'}\times[t_1,t_0]$, and let
$t_*\in(t_1,t_0]$ be its first vanishing time.
There would exist $z_*\in D'$ such that
\[
v(t_*,z_*)=0,
\qquad
v\ge0
\quad\text{on }\overline{D'}\times[t_1,t_*].
\]
Approaching $t_*$ from below gives
$\partial_tv(t_*,z_*)\le0$, while the interior
spatial minimum at $z_*$ gives
$(v_{i\bar j}(t_*,z_*))\ge0$.
Hence, at $(t_*,z_*)$,
\[
0\ge\partial_tv
\ge a^{i\bar j}v_{i\bar j}+\delta
\ge\delta>0,
\]
a contradiction.

Thus $v>0$ on
$\overline{D'}\times[t_1,t_0]$, so
\[
w(t_0,z)>-\delta(1+t_0)
\qquad(z\in D').
\]
The same inequality holds on $D\setminus D'$,
where $w(t_0,z)>-\delta$. Since $\delta>0$
and $t_0\in(0,\tau]$ are arbitrary, we conclude
that $U\ge V$ on $D\times(0,\tau]$.

The same proof, with $t_1=s_0$, shows that
$U\ge V$ on $D\times[s_0,t_0]$ whenever
$U(s_0)\ge V(s_0)$ on $D$ and
\eqref{eq:comparison-boundary} holds with
$I=[s_0,t_0]$.
\end{proof}

\begin{remark}
We call a smooth function with $t\theta_\mu+\ddbar U>0$ satisfying
\eqref{eq:comparison-super} a \emph{supersolution} of
\eqref{eq:scalar-gauge}, and use the term \emph{subsolution} for
\eqref{eq:comparison-sub}.
\end{remark}

\subsection{Potentials of positive closed currents}
\label{sec:currents}

We establish two consequences of
\eqref{eq:current-homogeneity}.
First, every current in
$\mathcal T_{Y,K}(\Omega)$ is exact and
hence admits a global plurisubharmonic
potential. Second, its Lelong numbers vanish.

\begin{proposition}\label{prop:current-exact-psh}
Let $T\in\mathcal T_{Y,K}(\Omega)$. Then
\begin{equation}\label{eq:current-cartan}
T=\frac1K d(\iota_YT).
\end{equation}
In particular, there exists
$\psi\in\PSH(\Omega)$ such that
\begin{equation}\label{eq:current-global-potential}
T=\ddbar\psi.
\end{equation}
The potential $\psi$ is unique up to the addition
of a real pluriharmonic function.
\end{proposition}

\begin{proof}
Differentiating \eqref{eq:current-homogeneity}
at $s=0$ in the sense of currents gives
$\Lie_YT=KT$. Since $dT=0$, Cartan's formula
yields
\[
KT=\Lie_YT
=d(\iota_YT)+\iota_Y(dT)
=d(\iota_YT).
\]
This proves \eqref{eq:current-cartan}. Equation~\eqref{eq:current-cartan} shows
that $T$ is exact as a current. By
\cite[Lemma~4.1]{Magnusson2011},
there exists a global plurisubharmonic
function $u$ such that $T=dd^c u$.
Since $dd^c=2\ddbar$, the function
$\psi=2u$ satisfies
\eqref{eq:current-global-potential}.
If $\psi_1$ and $\psi_2$ are two such
potentials, then
$\ddbar(\psi_1-\psi_2)=0$ in the distributional
sense. By Weyl's lemma
\cite[Theorem~IX.25]{ReedSimonII1975},
$\psi_1-\psi_2$ is a smooth real
pluriharmonic function.
\end{proof}

We next show that
\eqref{eq:current-homogeneity}
also implies the vanishing of Lelong numbers.

\begin{proposition}\label{prop:current-zero-lelong}
Let $T\in\mathcal T_{Y,K}(\Omega)$. Then
\[
\nu(T,p)=0
\qquad\text{for every }p\in\Omega.
\]
\end{proposition}

\begin{proof}
Lelong numbers are invariant under
biholomorphic pullback
\cite[Chapter~III, Corollary~7.2]{DemaillyCADG}.
Consequently,
\begin{equation}\label{eq:lelong-scaling}
\nu(T,\Phi_s(p))
=\nu(\Phi_s^*T,p)
=e^{Ks}\nu(T,p)
\qquad(s\in\R).
\end{equation}

Suppose that $a:=\nu(T,p)>0$ for some
$p\in\Omega$. Since $Y$ is real holomorphic,
the orbit map $s\mapsto\Phi_s(p)$ extends
to a holomorphic map
$\Psi:\Delta_\varepsilon\to\Omega$
for some $\varepsilon>0$, with
$\Psi(s)=\Phi_s(p)$ for real
$|s|<\varepsilon$. Choose $s_0\in(0,\varepsilon)$ and
$c$ satisfying
\[
a<c<e^{Ks_0}a.
\]
By \eqref{eq:lelong-scaling},
$\nu(T,\Psi(s))>c$ for all real $s$
in a sufficiently small interval
containing $s_0$.

By Siu's analyticity theorem
\cite{Siu1974Lelong}, the set
\[
E_c(T):=\{x\in\Omega:\nu(T,x)\ge c\}
\]
is a complex analytic subset of $\Omega$.
Thus $\Psi^{-1}(E_c(T))$ is an analytic
subset of $\Delta_\varepsilon$ containing
a nonempty real interval. Since every proper analytic subset of
a connected complex disc is discrete,
we obtain that $\Psi^{-1}(E_c(T))=\Delta_\varepsilon$.
In particular, $p=\Psi(0)\in E_c(T)$,
contradicting
$\nu(T,p)=a<c$.
\end{proof}

\subsection{Complete flows on exhaustion domains}
\label{sec:smooth-construction}

We first construct complete K\"ahler--Ricci flows
on smoothly bounded strongly pseudoconvex
subdomains of $\Omega$, 
starting from smooth semipositive
$(1,1)$-forms of the form $\ddbar\psi$.
The associated
scalar potentials tend to $+\infty$ at the
boundary, allowing us to compare the flows
on successive exhaustion domains. We will use
this construction in
Subsection~\ref{sec:rough-current}
to obtain flows starting from positive closed currents.

Let $\chi$ be a smooth strictly plurisubharmonic
exhaustion function on $\Omega$. Fix $z_0\in\Omega$
and choose regular values $c_j\uparrow\infty$
with $c_1>\chi(z_0)$. Let $\Omega_j$ be the
connected component of $\{\chi<c_j\}$
containing $z_0$. Since $\Omega$ is connected,
we have
\[
\Omega_j\Subset\Omega_{j+1},
\qquad
\bigcup_{j=1}^{\infty}\Omega_j=\Omega.
\]
Since $c_j$ is a regular value and $\chi$ is
strictly plurisubharmonic, each $\Omega_j$
has smooth strongly pseudoconvex boundary. Set $a=n+1$ and
\begin{equation}\label{eq:background-potential}
d_j=c_j-\chi,\qquad
\vartheta_j=-\log d_j+b_j,\qquad
\overline\omega_j=\ddbar\vartheta_j,
\end{equation}
where $b_j$ is a constant chosen so that $\vartheta_j\ge0$. We have
\begin{equation}\label{eq:logmetric-formula}
\overline\omega_j=
\frac{\ddbar\chi}{d_j}
+\frac{\ii\,\partial\chi\wedge\bar\partial\chi}{d_j^2}.
\end{equation}

Since
\[
|\partial\vartheta_j|_{\overline\omega_j}\le1
\]
and $\vartheta_j\to+\infty$ at
$\partial\Omega_j$, the metric
$\overline\omega_j$ is complete.

Moreover, $\overline\omega_j$ has bounded
geometry of every finite order. We argue as in
\cite[\S1, Proposition~1.3]{CY80}.
Fix $j$. Near $\partial\Omega_j$, choose
holomorphic coordinates centered at
$p\in\Omega_j$ such that $\partial\chi(p)$
has only a complex normal component.
Put $d=d_j(p)$ and rescale by
\[
(z',z_n)=(\sqrt d\,w',dw_n).
\]
Writing $D(w)=d_j(z(w))/d$, we have
$D(0)=1$, the tangential first derivatives
of $D$ vanish at the origin, and its complex
normal derivative is bounded away from zero
for $p$ sufficiently close to the boundary.
The smoothness of $\chi$ near
$\overline\Omega_j$ gives uniform bounds
for all derivatives of $D$ of every fixed
finite order on a sufficiently small ball.
This ball can be chosen with radius independent
of $p$, so that $D$ remains uniformly positive
and its complex normal derivative remains
bounded away from zero.

In these coordinates, the first term of
\eqref{eq:logmetric-formula} controls the
tangential directions, while the second
term controls the complex normal direction.
The strict plurisubharmonicity of $\chi$
therefore gives uniform upper and lower
bounds for the coefficients of
$\overline\omega_j$, together with their
derivatives of every finite order, on
these fixed-radius balls. The remaining
compact subset of $\Omega_j$ is covered
by finitely many ordinary coordinate charts.
Smooth functions and covariant tensors defined
near $\overline\Omega_j$ have uniformly
bounded coordinate derivatives in these charts.
All constants and chart radii may depend
on $j$, but not on the chart center.

We also have
\begin{equation}\label{eq:baromega-det}
\overline\omega_j^n
=e^{a\vartheta_j+\ell_j}\mu,
\qquad
\ell_j\in C^\infty(\overline\Omega_j).
\end{equation}
Indeed, writing $A=(\chi_{i\bar k})$
and $v=(\chi_i)$ in local holomorphic
coordinates, the matrix of
$\overline\omega_j$ is
$d_j^{-1}A+d_j^{-2}vv^*$, and hence
\[
\begin{aligned}
\det(d_j^{-1}A+d_j^{-2}vv^*)
&=d_j^{-n}\det A
\left(1+d_j^{-1}v^*A^{-1}v\right) =d_j^{-a}
\left(d_j\det A+
v^*\operatorname{adj}(A)v\right).
\end{aligned}
\]
The factor in parentheses extends
smoothly to $\overline\Omega_j$ and
is strictly positive in $\Omega_j$.
It is also positive on
$\partial\Omega_j$, since $v\ne0$
there and $\operatorname{adj}(A)>0$.
As $\mu$ is a smooth positive volume
form near $\overline\Omega_j$, this
proves \eqref{eq:baromega-det}.
In particular, $\ell_j$ is bounded
on $\Omega_j$. Taking the Ricci form
gives
\begin{equation}\label{eq:background-ricci}
\Ric_{\overline\omega_j}
=-\theta-a\overline\omega_j-\ddbar\ell_j.
\end{equation}

\begin{proposition}\label{prop:fixed-domain-flow}
Let $\psi\in C^\infty(\Omega)\cap\PSH(\Omega)$,
and set $\eta=\ddbar\psi$.
For every $j$, there exists a K\"ahler--Ricci
flow $\omega_j(t)$ on $\Omega_j\times(0,\infty)$
with the following properties:
\begin{enumerate}[label=\textup{(\roman*)}]
\item $\omega_j(t)$ is complete for every $t>0$.

\item There is a smooth function $U_j$ on
$\Omega_j\times(0,\infty)$ satisfying
\[
\omega_j(t)=t\theta+\ddbar U_j(t),
\qquad
\partial_tU_j
=\log\frac{\omega_j(t)^n}{\mu}.
\]

\item $U_j(t)\to\psi$ locally uniformly on
$\Omega_j$ as $t\downarrow0$.

\item There exists $C_j>0$ such that
\begin{equation}\label{eq:fixed-corner}
U_j(t,z)\ge
\psi(z)+at\vartheta_j(z)
+nt\log(at)-nt-C_jt
\end{equation}
for every $z\in\Omega_j$ and $t>0$.
Consequently,
\begin{equation}\label{eq:fixed-boundary-blow}
U_j(t,z)\longrightarrow+\infty
\qquad(z\to\partial\Omega_j),
\end{equation}
locally uniformly for $t\in(0,\infty)$.
\end{enumerate}

Moreover, there exists $B_j\ge1$ such that,
for every $0<\varepsilon<1$, there is a
complete K\"ahler--Ricci flow
$\omega_{j,\varepsilon}(t)$ on
$\Omega_j\times[0,\infty)$, with bounded
curvature on each compact time interval and
initial metric
\begin{equation}\label{eq:fixed-regularized-initial}
\omega_{j,\varepsilon}(0)
=\eta+\varepsilon B_j\overline\omega_j.
\end{equation}
For some sequence $\varepsilon_k\downarrow0$,
\[
\omega_{j,\varepsilon_k}\longrightarrow\omega_j
\quad\text{in }
C^\infty_{\mathrm{loc}}
(\Omega_j\times(0,\infty)).
\]
For every compact interval $I\subset(0,\infty)$,
there exists $C_{j,I}>0$, independent of
$\varepsilon$, such that
\begin{equation}\label{eq:fixed-regularized-equivalence}
C_{j,I}^{-1}\overline\omega_j
\le\omega_{j,\varepsilon}(t)
\le C_{j,I}\overline\omega_j,
\qquad t\in I.
\end{equation}
\end{proposition}

\begin{proof}
Fix $j$. Since $\ddbar\chi$ is positive definite and
$\eta$ is smooth and semipositive on
$\overline\Omega_j$, there exists $C_j>0$
such that
\[
0\le\eta\le C_j\ddbar\chi.
\]
As $d_j$ is bounded above on $\Omega_j$,
\eqref{eq:logmetric-formula} allows us to
choose $B_j\ge1$ such that
\begin{equation}\label{eq:fixed-initial-upper}
0\le\eta\le h_j,\qquad h_j=B_j\overline\omega_j.
\end{equation}

The metric $h_j$ is complete and has bounded
geometry of every finite order. For
$0<\varepsilon<1$, the metric
$\eta+\varepsilon h_j$ satisfies
\[
\varepsilon h_j\le\eta+\varepsilon h_j
\le(1+\varepsilon)h_j.
\]
Since $\eta$ is smooth on a neighborhood of
$\overline\Omega_j$, the bounded-geometry charts
also give bounded curvature for this metric. For every $S>0$, \eqref{eq:background-ricci} gives
\begin{equation}\label{eq:fixed-chaulee-condition}
\begin{aligned}
\eta-S\Ric_{h_j}+S\ddbar(-\ell_j)
&=\eta+S\theta+aS\overline\omega_j \ge\frac{aS}{B_j}h_j
>\frac{aS}{2B_j}h_j.
\end{aligned}
\end{equation}
Since $\ell_j$ is bounded, we may apply
\cite[Theorem~1.4]{ChauLee2020} with
$f=-\ell_j$.

Fix $S>0$ and $0<\varepsilon<1$.
By \cite[Theorem~1.1]{ChauLee2020},
the complete bounded-curvature metric
$\eta+\varepsilon h_j$ admits a complete
bounded-curvature K\"ahler--Ricci flow on
a short time interval. Let
$[0,T_{j,\varepsilon})$ be its maximal
interval of existence in this class.

Define
\begin{equation}\label{eq:fixed-chaulee-potential}
\phi_{j,\varepsilon}(t)
=\int_0^t
\log\frac{\omega_{j,\varepsilon}(s)^n}{h_j^n}\,ds.
\end{equation}
Let $0<b<\min\{S,T_{j,\varepsilon}\}$.
On $[0,b]$, the flow has bounded curvature
and is uniformly equivalent to $h_j$.
By \cite[\S7, Claim~7.1]{ChauLee2020},
\[
nt\log t-C_{j,S}t
\le\phi_{j,\varepsilon}(t)
\le C_{j,S}t,
\qquad
\partial_t\phi_{j,\varepsilon}(t)
\le C_{j,S},
\]
where $C_{j,S}$ is independent of
$\varepsilon$ and $b$. Write $\omega=\omega_{j,\varepsilon}$,
$\phi=\phi_{j,\varepsilon}$,
$H=\tr_\omega h_j$, and
$P=\partial_t-\Delta_\omega$.
Set
\[
V=(S-t)\partial_t\phi+\phi+S\ell_j+nt.
\]
By \eqref{eq:KRF},
\eqref{eq:fixed-regularized-initial}, and
\eqref{eq:fixed-chaulee-potential},
\[
\omega=\eta+\varepsilon h_j
-t\Ric_{h_j}+\ddbar\phi.
\]
Consequently,
\[
\begin{aligned}
P(\partial_t\phi)
&=-\tr_\omega\Ric_{h_j},\\
P\phi
&=\partial_t\phi-n+
\tr_\omega(\eta+\varepsilon h_j-t\Ric_{h_j}).
\end{aligned}
\]
Using \eqref{eq:fixed-chaulee-condition},
we obtain
\[
PV
=\tr_\omega\bigl(
\eta+\varepsilon h_j-S\Ric_{h_j}
-S\ddbar\ell_j\bigr)
\ge\beta_S H,
\qquad
\beta_S=\frac{aS}{2B_j}.
\]
Moreover, the parabolic trace inequality
\cite[(2.2)]{ChauLee2020} gives
\[
P\log H\le C_{h_j}H,
\]
where $C_{h_j}$ depends only on the
curvature bound of $h_j$. Choose $L>0$ with
$L\beta_S\ge C_{h_j}+2$, and put
$m=1+LnS$. For $0<\delta\le1$, set
\[
G_\delta
=\log H-LV+m\log t-\delta\vartheta_j.
\]
Since
$P(-\delta\vartheta_j)=\delta H/B_j$,
we have
\[
PG_\delta\le \frac mt-H.
\]
For fixed $\varepsilon$ and $b$,
the functions $H$, $\phi$, and
$\partial_t\phi$ are bounded on
$\Omega_j\times[0,b]$.
As $\vartheta_j\to+\infty$ at
$\partial\Omega_j$, the function
$G_\delta$ attains its maximum on
$\Omega_j\times(0,b]$.
At a maximum point, $tH\le m$. 

The arithmetic--geometric mean inequality
gives
\[
-\partial_t\phi\le n\log(H/n).
\]
Using the lower bound for $\phi$,
the boundedness of $\ell_j$, and
$\vartheta_j\ge0$, we obtain
at the maximum point
\[
\begin{aligned}
G_\delta
&\le [1+Ln(S-t)]\log H
+(m-Lnt)\log t+C_{j,S}\\
&=(m-Lnt)\log(tH)+C_{j,S}\\
&\le m\log m+C_{j,S}.
\end{aligned}
\]
Here the last inequality follows from
$tH\le m$ and $1\le m-Lnt\le m$. The constants are independent of
$\varepsilon$, $b$, and $\delta$.
Letting $\delta\downarrow0$ and using
the upper bounds for $\phi$ and
$\partial_t\phi$, we obtain
\[
H\le C_{j,S}t^{-m}
\qquad (0<t\le b).
\]
Since $b$ is arbitrary, for every
$t_0>0$ we conclude that
\[
\tr_{\omega_{j,\varepsilon}(t)}h_j
\le C_{j,S,t_0},
\qquad
t_0\le t<\min\{S,T_{j,\varepsilon}\},
\]
with $C_{j,S,t_0}$ independent
of $\varepsilon$. In particular,
\[
\frac{\omega_{j,\varepsilon}(t)^n}{h_j^n}
\le e^{C_{j,S}}.
\]

Let $\lambda_1,\ldots,\lambda_n$ denote
the eigenvalues of $\omega_{j,\varepsilon}(t)$
relative to $h_j$. Then
\[
\sum_{p=1}^n\lambda_p^{-1}
\le C_{j,S,t_0},
\qquad
\prod_{p=1}^n\lambda_p\le e^{C_{j,S}}.
\]
Consequently,
\[
C_{j,S,t_0}^{-1}h_j
\le\omega_{j,\varepsilon}(t)
\le e^{C_{j,S}}C_{j,S,t_0}^{n-1}h_j
\]
throughout the same time interval.

Suppose that $T:=T_{j,\varepsilon}<S$.
Taking $t_0=T/2$ in the preceding
estimates, we obtain $C\ge1$ such that
\[
C^{-1}h_j\le\omega_{j,\varepsilon}(t)
\le Ch_j,
\qquad T/2\le t<T.
\]
Set $\widehat h_j=C^{-1}h_j$.
Since $h_j$ is complete and has bounded
curvature, so does $\widehat h_j$, and
\[
\widehat h_j\le\omega_{j,\varepsilon}(t)
\le C^2\widehat h_j,
\qquad T/2\le t<T.
\]
By \cite[Theorem~1.1]{ChauLee2020},
there exists $\tau>0$, depending only
on $n$, $C$, and the curvature bound
of $\widehat h_j$, such that for every
$t_*\in[T/2,T)$, the K\"ahler--Ricci
flow starting from
$\omega_{j,\varepsilon}(t_*)$
admits a complete bounded-curvature
solution on $[t_*,t_*+\tau)$.
Choose $t_*$ sufficiently close to
$T$ that $t_*+\tau>T$.
By \cite[Theorem~1.1]{ChenZhuRicciFlowUniqueness},
with Riemannian Ricci-flow time $s=t/2$,
this solution agrees with
$\omega_{j,\varepsilon}$ on
$[t_*,T)$. It therefore extends
$\omega_{j,\varepsilon}$ past $T$,
contrary to maximality.
Hence $T_{j,\varepsilon}\ge S$.
Since $S>0$ was arbitrary,
$T_{j,\varepsilon}=\infty$.
For every compact interval
$I\subset(0,\infty)$, the preceding
two-sided bounds, together with
$h_j=B_j\overline\omega_j$, yield
\eqref{eq:fixed-regularized-equivalence},
with a constant independent of
$\varepsilon$.

Let $E\Subset\Omega_j$. Applying
\cite[Corollary~1.2]{ShermanWeinkove2012}
on bounded-geometry coordinate balls
to the time-shifted flows
$\omega_{j,\varepsilon}(t_0/2+s)$,
and using \eqref{eq:fixed-regularized-equivalence}
on $[t_0/2,t_1]$, gives uniform bounds
for their spatial derivatives of every
order on $E\times[t_0,t_1]$.
Differentiating the K\"ahler--Ricci flow equation
\eqref{eq:KRF} then gives the corresponding
mixed space--time derivative bounds.
By Arzel\`a--Ascoli and a diagonal
subsequence argument, there exists
$\varepsilon_k\downarrow0$ such that
\begin{equation}\label{eq:fixed-flow-smooth-limit}
\omega_{j,\varepsilon_k}\longrightarrow\omega_j
\quad\text{in }C^\infty_{\mathrm{loc}}
(\Omega_j\times(0,\infty)).
\end{equation}
The limit satisfies \eqref{eq:KRF}, and
the estimates \eqref{eq:fixed-regularized-equivalence}
pass to the limit. Thus $\omega_j(t)$
is a complete K\"ahler metric for every $t>0$. We now construct the scalar potential.
For $0<\varepsilon<1$, set
\[
v_{j,\varepsilon}(t,z)
=\int_0^t
\left(
\log\frac{\omega_{j,\varepsilon}(s,z)^n}
{\overline\omega_j(z)^n}
+\ell_j(z)
\right)\,ds.
\]

Since
\[
\Ric_{\omega_{j,\varepsilon}}
=\Ric_{\overline\omega_j}
-\ddbar\log\frac{\omega_{j,\varepsilon}^n}
{\overline\omega_j^n},
\]
equation \eqref{eq:background-ricci} gives
\[
-\Ric_{\omega_{j,\varepsilon}}
=\theta+a\overline\omega_j+
\ddbar\left(
\log\frac{\omega_{j,\varepsilon}^n}
{\overline\omega_j^n}+\ell_j
\right).
\]
Integrating \eqref{eq:KRF} from $0$ to $t$ gives
\begin{align}
\omega_{j,\varepsilon}(t)
&=\eta+t\theta+
(B_j\varepsilon+at)\overline\omega_j
+\ddbar v_{j,\varepsilon}(t),
\label{eq:fixed-regularized-form}\\
\partial_t v_{j,\varepsilon}
&=\log\frac{\omega_{j,\varepsilon}^n}
{\overline\omega_j^n}+\ell_j,
\qquad v_{j,\varepsilon}(0)=0.
\label{eq:fixed-regularized-scalar}
\end{align}

Choose $M_j>0$ such that
$\theta\le M_j\overline\omega_j$
on $\Omega_j$. This is possible since
$\theta$ is smooth near
$\overline\Omega_j$ and
$\overline\omega_j$ dominates a positive
multiple of $\ddbar\chi$.
Set
\[
\begin{aligned}
q^-_{j,\varepsilon}(t)
&=n\int_0^t
\log(B_j\varepsilon+as)\,ds
+t\inf_{\Omega_j}\ell_j,\\
q^+_{j,\varepsilon}(t)
&=n\int_0^t
\log\bigl(B_j(1+\varepsilon)+(a+M_j)s\bigr)\,ds
+t\sup_{\Omega_j}\ell_j.
\end{aligned}
\]
By \eqref{eq:fixed-initial-upper} and
$0\le\theta\le M_j\overline\omega_j$,
these functions are respectively a
subsolution and a supersolution of
\eqref{eq:fixed-regularized-scalar},
both with zero initial value.

Fix $\varepsilon\in(0,1)$ and $T>0$.
Since $\omega_{j,\varepsilon}$ is a complete
bounded-curvature flow on $[0,T]$, it is
uniformly equivalent to $\overline\omega_j$
there, with constants depending on
$j,\varepsilon,T$. In particular,
$v_{j,\varepsilon}$ is bounded on
$\Omega_j\times[0,T]$.

For either
\[
w=v_{j,\varepsilon}-q^+_{j,\varepsilon}
\quad\text{or}\quad
w=q^-_{j,\varepsilon}-v_{j,\varepsilon},
\]
linearizing the logarithmic determinant
gives
\[
\partial_tw\le a^{p\bar q}w_{p\bar q},
\]
where $(a^{p\bar q})$ is positive
definite and uniformly bounded
relative to $\overline\omega_j^{-1}$
on $[0,T]$, for fixed $\varepsilon$.
Consequently,
\[
a^{p\bar q}(\vartheta_j)_{p\bar q}
\le A_{j,\varepsilon,T}
\]
for some constant $A_{j,\varepsilon,T}>0$. For $\delta>0$, set
\[
w_\delta=w-\delta\vartheta_j
-\delta(A_{j,\varepsilon,T}+1)t.
\]
Then
\[
(\partial_t-a^{p\bar q}\partial_p\partial_{\bar q})
w_\delta\le-\delta.
\]

Since $w$ is bounded on
$\Omega_j\times[0,T]$ and
$\vartheta_j\to+\infty$ at
$\partial\Omega_j$, the function
$w_\delta$ tends uniformly to $-\infty$
at the boundary. Thus, for sufficiently
large regular values $R$, we have
$w_\delta<0$ on
$\{\vartheta_j=R\}\times[0,T]$.
Since $w_\delta(0,\cdot)\le0$ and
\[
(\partial_t-a^{p\bar q}
\partial_p\partial_{\bar q})w_\delta
\le-\delta,
\]
the weak parabolic maximum principle
\cite[\S7.1.4, Theorem~8]{PDEbookEvans_1998},
applied on $\{\vartheta_j<R\}$,
gives $w_\delta\le0$ there.
Letting $R\to\infty$ and then
$\delta\downarrow0$ yields $w\le0$.
See also \cite[\S7, proof of
Claim~7.1]{ChauLee2020} for the
same exhaustion argument. As $T$ was
arbitrary, we obtain
\begin{equation}\label{eq:fixed-regularized-v-bounds}
q^-_{j,\varepsilon}(t)
\le v_{j,\varepsilon}(t,z)
\le q^+_{j,\varepsilon}(t)
\qquad(t\ge0,\ z\in\Omega_j).
\end{equation}
Since $0<\varepsilon<1$,
\[
q^-_{j,\varepsilon}(t)
\ge nt\log(at)-nt+
t\inf_{\Omega_j}\ell_j,
\]
whereas $q^+_{j,\varepsilon}(t)\le C_jt$
for $0<t\le1$, with $C_j$ independent
of $\varepsilon$. Therefore, after
enlarging $C_j$,
\begin{equation}\label{eq:fixed-uniform-small-time}
\sup_{\Omega_j}|v_{j,\varepsilon}(t)|
\le C_jt(1+|\log t|)
\qquad(0<t\le1).
\end{equation}

We pass to the sequence
$\varepsilon_k\downarrow0$ in
\eqref{eq:fixed-flow-smooth-limit}.
By \eqref{eq:fixed-regularized-scalar},
\[
\partial_t v_{j,\varepsilon_k}
\longrightarrow
\log\frac{\omega_j^n}{\overline\omega_j^n}
+\ell_j
\quad\text{in }C^\infty_{\mathrm{loc}}
(\Omega_j\times(0,\infty)).
\]
Fix $E\Subset\Omega_j$ and
$0<t_0<t_1$. For
$0<\delta<\min\{1,t_0\}$,
integrating the difference of the
time derivatives from $\delta$ to
$t$ and using
\eqref{eq:fixed-uniform-small-time}
gives
\[
\begin{aligned}
&\sup_{E\times[t_0,t_1]}
|v_{j,\varepsilon_k}-v_{j,\varepsilon_l}| \le
2C_j\delta(1+|\log\delta|) +
(t_1-\delta)
\sup_{E\times[\delta,t_1]}
|\partial_tv_{j,\varepsilon_k}
-\partial_tv_{j,\varepsilon_l}|.
\end{aligned}
\]
Letting first $k,l\to\infty$ and
then $\delta\downarrow0$, we see
that $v_{j,\varepsilon_k}$ converges
locally uniformly to a function $v_j$.

Taking the trace of
\eqref{eq:fixed-regularized-form}
with respect to $\overline\omega_j$
and subtracting the identities for
$\varepsilon_k$ and $\varepsilon_l$,
we obtain
\[
\begin{aligned}
&\Delta_{\overline\omega_j}
(v_{j,\varepsilon_k}-v_{j,\varepsilon_l}) \quad=
\tr_{\overline\omega_j}
(\omega_{j,\varepsilon_k}
-\omega_{j,\varepsilon_l})
-nB_j(\varepsilon_k-\varepsilon_l).
\end{aligned}
\]

The right-hand side tends to zero in
$C^\infty_{\mathrm{loc}}
(\Omega_j\times(0,\infty))$,
while $v_{j,\varepsilon_k}-v_{j,\varepsilon_l}$
tends to zero locally uniformly.
Since $\Delta_{\overline\omega_j}$ is a fixed
elliptic operator with smooth coefficients,
the interior estimates
\cite[\S6.3.1, Theorem~2]{PDEbookEvans_1998},
applied in local coordinates, and Sobolev
embedding imply convergence of all spatial
derivatives on compact subsets. The constants
are uniform on compact time intervals,
since $\Delta_{\overline\omega_j}$
is independent of time.
By \eqref{eq:fixed-regularized-scalar}
and \eqref{eq:fixed-flow-smooth-limit},
all mixed space--time derivatives
also converge. Hence
\[
v_{j,\varepsilon_k}\longrightarrow v_j
\quad\text{in }C^\infty_{\mathrm{loc}}
(\Omega_j\times(0,\infty)).
\]

Passing to the limit in
\eqref{eq:fixed-regularized-form}
and \eqref{eq:fixed-regularized-scalar}
gives
\begin{align}
\omega_j(t)
&=\eta+t\theta+
at\overline\omega_j+\ddbar v_j(t),
\label{eq:fixed-flow-potential}\\
\partial_tv_j
&=\log\frac{\omega_j^n}
{\overline\omega_j^n}+\ell_j.
\label{eq:fixed-flow-scalar}
\end{align}
Letting $\varepsilon_k\downarrow0$
in \eqref{eq:fixed-regularized-v-bounds},
we obtain
\begin{equation}\label{eq:fixed-potential-bounds}
\begin{aligned}
nt\log(at)-nt+t\inf_{\Omega_j}\ell_j
&\le v_j(t,z) \le n\int_0^t
\log\bigl(B_j+(a+M_j)s\bigr)\,ds
+t\sup_{\Omega_j}\ell_j.
\end{aligned}
\end{equation}
Both bounds tend to zero uniformly
in $z\in\Omega_j$ as $t\downarrow0$.
Thus
\begin{equation}\label{eq:fixed-v-initial}
\lim_{t\downarrow0}
\sup_{\Omega_j}|v_j(t)|=0.
\end{equation}

Define
\[
U_j(t)=\psi+at\vartheta_j+v_j(t).
\]
Since $\eta=\ddbar\psi$ and
$\overline\omega_j=\ddbar\vartheta_j$,
\eqref{eq:fixed-flow-potential} gives
\[
\omega_j(t)=t\theta+\ddbar U_j(t).
\]
Moreover, by \eqref{eq:baromega-det}
and \eqref{eq:fixed-flow-scalar},
\[
\begin{aligned}
\partial_tU_j
&=a\vartheta_j+
\log\frac{\omega_j^n}{\overline\omega_j^n}
+\ell_j =\log\frac{\omega_j^n}{\mu}.
\end{aligned}
\]
Thus $U_j$ satisfies
\eqref{eq:scalar-gauge}.
Since $\vartheta_j$ is bounded on
compact subsets of $\Omega_j$,
\eqref{eq:fixed-v-initial} gives
$U_j(t)\to\psi$ locally uniformly
as $t\downarrow0$.

The lower bound in
\eqref{eq:fixed-potential-bounds}
gives \eqref{eq:fixed-corner}
after choosing
$C_j\ge\max\{1,-\inf_{\Omega_j}\ell_j\}$.
Finally, $\psi$ is bounded on
$\overline\Omega_j$, and
$\vartheta_j\to+\infty$ at
$\partial\Omega_j$. The remaining
terms in \eqref{eq:fixed-corner}
are bounded below on every compact
positive-time interval. This proves
\eqref{eq:fixed-boundary-blow}.
\end{proof}

We shall also use the complete K\"ahler--Einstein
metrics on the exhaustion domains $\Omega_j$.
The following lemma gives their comparison
with $\overline\omega_j$ and their convergence
to the fixed metric $\omega_E$ on $\Omega$.

\begin{lemma}\label{lem:stein-KE-exhaustion}
Let $\omega_E$ be the complete K\"ahler--Einstein
metric on $\Omega$ fixed above, satisfying
$\Ric_{\omega_E}=-K\omega_E$.
For every $j$, there exists a complete
K\"ahler--Einstein metric $\eta_j$ on $\Omega_j$
such that
\[
\Ric_{\eta_j}=-K\eta_j,
\qquad
C_j^{-1}\overline\omega_j
\le\eta_j\le C_j\overline\omega_j
\]
for some $C_j\ge1$. Moreover,
\[
\omega_E^n\le\eta_{j+1}^n\le\eta_j^n
\qquad\text{on }\Omega_j.
\]
For every $D\Subset\Omega$, there exist
$j_D$ and $C_D>0$ such that
\[
C_D^{-1}h\le\eta_j\le C_Dh
\qquad\text{on }D,\quad j\ge j_D.
\]
The metrics $\eta_j$ converge to $\omega_E$
in $C^\infty_{\mathrm{loc}}(\Omega)$.
\end{lemma}

\begin{proof}
Fix $j$ and set
\[
\alpha_j=\overline\omega_j+\frac{\theta}{a},
\qquad a=n+1.
\]
Since $\theta$ is smooth near
$\overline\Omega_j$, and
$\overline\omega_j$ dominates a positive
multiple of $\ddbar\chi$, there exists
$C_j>0$ such that
\[
\overline\omega_j\le\alpha_j
\le C_j\overline\omega_j.
\]
In particular, $\alpha_j$ is complete.
The bounded-geometry charts for
$\overline\omega_j$ constructed after
\eqref{eq:logmetric-formula}, together
with the smoothness of $\theta$ near
$\overline\Omega_j$, show that
$\alpha_j$ has bounded geometry
of every finite order.

We claim that
\[
\alpha_j^n=e^{a\vartheta_j+\ell'_j}\mu,
\qquad
\ell'_j\in C^\infty(\overline\Omega_j).
\]
Indeed, in local holomorphic coordinates,
put
\[
A_j=\left(\chi_{p\bar q}
+\frac{d_j}{a}\theta_{p\bar q}\right),
\qquad
v=(\chi_p).
\]
Then
\[
\det\alpha_j
=d_j^{-a}
\left(d_j\det A_j+
v^*\operatorname{adj}(A_j)v\right).
\]
The factor in parentheses is smooth and
strictly positive on $\overline\Omega_j$,
since $A_j$ is positive definite and
$d\chi\ne0$ on $\partial\Omega_j$.
This proves the asserted volume identity.

Since $\ell'_j$ extends smoothly to
$\overline\Omega_j$, the coordinate
estimates established after
\eqref{eq:logmetric-formula} give uniform
$C^{k-2,\gamma}$ bounds in the
bounded-geometry charts, for every
$k\ge5$ and $0<\gamma<1$.
Thus \cite[Theorem~4.4]{CY80} applies
to the complete metric $\alpha_j$,
with coefficient $a>0$ and right-hand
side $-\ell'_j$. It gives a solution
$f_j$ satisfying
\[
(\alpha_j+\ddbar f_j)^n
=e^{af_j-\ell'_j}\alpha_j^n,
\qquad
C_j^{-1}\alpha_j
\le\alpha_j+\ddbar f_j
\le C_j\alpha_j.
\]
Taking $k$ arbitrarily large and using
uniqueness in \cite[Theorem~4.4]{CY80} gives $f_j\in C^\infty(\Omega_j)$.

Since $\theta=-\Ric_h$ and
$\alpha_j=\overline\omega_j+\theta/a$,
the volume identity yields
\[
\Ric_{\alpha_j}
=-a\alpha_j-\ddbar\ell'_j.
\]
Consequently,
\[
\Ric_{\alpha_j+\ddbar f_j}
=-a(\alpha_j+\ddbar f_j).
\]
It follows that
\[
\eta_j=\frac aK(\alpha_j+\ddbar f_j)
\]
is a complete K\"ahler--Einstein metric
satisfying $\Ric_{\eta_j}=-K\eta_j$.
Combining the bounds for $\alpha_j$
and $\alpha_j+\ddbar f_j$, and enlarging
$C_j$ if necessary, we obtain
\[
C_j^{-1}\overline\omega_j
\le\eta_j\le C_j\overline\omega_j.
\]

We next use the Schwarz lemma for volume
forms \cite[\S1.1]{MokYau}.
If $\gamma$ and $\tau$ are K\"ahler
metrics satisfying
\[
\Ric_\gamma=-K\gamma,
\qquad
\Ric_\tau=-K\tau,
\]
and $\gamma$ is complete, then
\[
\tau^n\le\gamma^n.
\]
Indeed, the scalar curvature
of $\gamma$ is $-nK$, while
\[
\ddbar\log\tau^n=K\tau,
\qquad
(\ddbar\log\tau^n)^n=K^n\tau^n.
\]
Applying the volume Schwarz lemma with
$K_1=nK$ and any $0<K_2<K^n$ gives
\[
\frac{\tau^n}{\gamma^n}
\le\frac{K^n}{K_2}.
\]
Letting $K_2\uparrow K^n$ proves
the claimed inequality.

Applying this comparison on $\Omega_j$ with
$(\gamma,\tau)=(\eta_j,\eta_{j+1})$, and on
$\Omega_{j+1}$ with
$(\gamma,\tau)=(\eta_{j+1},\omega_E)$,
we obtain
\[
\omega_E^n\le\eta_{j+1}^n
\le\eta_j^n
\qquad\text{on }\Omega_j.
\]

Set $a_n=(n+1)\kappa/(2n)$.
Since $\Ric_{\eta_j}=-K\eta_j$,
the family $(1+Kt)\eta_j$, $t\ge0$,
is a complete K\"ahler--Ricci flow.
The holomorphic sectional curvature
of $h$ is at most $-\kappa$.
Applying \cite[Lemma~2.1]{L22}
with reference metric $h$, we obtain
\[
(1+Kt)\eta_j\ge a_nt\,h.
\]
Dividing by $1+Kt$ and letting
$t\to\infty$ gives
\[
\eta_j\ge\frac{a_n}{K}h
\qquad\text{on }\Omega_j.
\]

Fix $D\Subset\Omega$, and choose
$j_D$ such that
$\overline D\Subset\Omega_{j_D}$.
By the volume comparison,
\[
\frac{\eta_j^n}{h^n}
\le\frac{\eta_{j_D}^n}{h^n}
\le M_D
\qquad\text{on }D,\quad j\ge j_D,
\]
where
\[
M_D=\sup_D\frac{\eta_{j_D}^n}{h^n}<\infty.
\]
If $\lambda_1,\ldots,\lambda_n$
are the eigenvalues of $\eta_j$
relative to $h$, then
\[
\lambda_p\ge\frac{a_n}{K},
\qquad
\prod_{p=1}^n\lambda_p\le M_D.
\]
Hence
\[
\frac{a_n}{K}h
\le\eta_j
\le M_D\left(\frac K{a_n}\right)^{n-1}h
\qquad\text{on }D,\quad j\ge j_D.
\]

Fix $D\Subset D'\Subset\Omega$.
The preceding estimates imply that
the complete flows
$(1+Kt)\eta_j$ are uniformly equivalent
to $h$ on $D'\times[1/2,2]$,
for all sufficiently large $j$.

The interior derivative estimates
\cite[Theorem~1.1 and Corollary~1.2]{ShermanWeinkove2012},
applied in fixed holomorphic coordinate balls
contained in $D'$, give uniform bounds
for the spatial derivatives of $\eta_j$
of every order on $D$.
Thus $\{\eta_j\}$ is precompact in
$C^\infty_{\mathrm{loc}}(\Omega)$.

Let $\eta_\infty$ be a subsequential
limit. The two-sided metric bounds
show that $\eta_\infty$ is positive
definite, and smooth convergence gives
\[
\Ric_{\eta_\infty}=-K\eta_\infty.
\]
The volume comparison implies
\[
\omega_E^n\le\eta_\infty^n.
\]
Applying the same comparison with
complete source $\omega_E$ and
target $\eta_\infty$ gives the reverse
inequality. Therefore
$\eta_\infty^n=\omega_E^n$, and
\[
-K\eta_\infty
=-\ddbar\log\eta_\infty^n
=-\ddbar\log\omega_E^n
=-K\omega_E.
\]
Hence $\eta_\infty=\omega_E$.
Since every convergent subsequence
has the same limit, the entire
sequence converges to $\omega_E$
in $C^\infty_{\mathrm{loc}}(\Omega)$.
\end{proof}

\subsection{K\"ahler--Ricci flow from a positive current}
\label{sec:rough-current}

Let $T=\ddbar\psi\ge0$, where $\psi\in\PSH(\Omega)$ and $\nu(T,z)=0$ for every $z\in\Omega$. We construct a K\"ahler--Ricci flow starting from $T$ by applying Proposition~\ref{prop:fixed-domain-flow} to decreasing smooth approximations of $\psi$ and passing to the exhaustion limit.

Let $\omega_E$ be the metric of Lemma~\ref{lem:stein-KE-exhaustion}, and set
\begin{equation}\label{eq:KE-potential-rho}
\rho=\frac1K\log\frac{\omega_E^n}{\mu}.
\end{equation}
The Ricci identity gives $\theta+K\ddbar\rho=K\omega_E$.
We use the following local estimates.

\begin{lemma}\label{lem:W2p-L1}
Let $B_r\Subset B_R$ be concentric balls in $\R^N$, and let $1<p<\infty$. There exists $C=C(r,R,p,N)>0$ such that
\[
\|u\|_{W^{2,p}(B_r)}
\le C\left(\|\Delta u\|_{L^p(B_R)}
+\|u\|_{L^1(B_R)}\right)
\]
for every $u\in L^1(B_R)$ whose distributional Laplacian belongs to $L^p(B_R)$.
\end{lemma}

\begin{proof}
Choose $r<r_1<R$ and $\zeta\in C_c^\infty(B_R)$
with $\zeta=1$ on $B_{r_1}$. Extend
$f=\zeta\Delta u$ by zero to $\R^N$ and let $v$ be its
Newtonian potential \cite[(9.26)]{GilbargTrudinger2001}.
By \cite[Theorem~9.9]{GilbargTrudinger2001} for $N\ge2$,
and by direct integration when $N=1$, $\Delta v=f$ almost everywhere and
\[
\|D^2v\|_{L^p(B_{r_1})}
\le C\|\Delta u\|_{L^p(B_R)}.
\]
Since $f$ has compact support and the
fundamental solution $\Gamma$ of the Laplacian,
given by \cite[(2.12)]{GilbargTrudinger2001} for $N\ge2$ and by
$\Gamma(x)=|x|/2$ for $N=1$,
satisfies $\Gamma,D\Gamma\in L^1_{\mathrm{loc}}
(\R^N)$, Young's convolution inequality gives
\[
\|v\|_{W^{1,p}(B_{r_1})}
\le C\|\Delta u\|_{L^p(B_R)}.
\]
Combining the last two estimates of $v$ gives 
\[
\|v\|_{W^{2,p}(B_{r_1})}
\le C\|\Delta u\|_{L^p(B_R)}.
\]

Since $\Delta(u-v)=0$ distributionally on
$B_{r_1}$, Weyl's lemma
\cite[Theorem~IX.25]{ReedSimonII1975}
implies that $w=u-v$ is harmonic there.
Choose $r<s<r_1$. The mean-value property
\cite[Theorem~2.1]{GilbargTrudinger2001}
gives
\[
\sup_{B_s}|w|
\le C\|w\|_{L^1(B_{r_1})}.
\]
The interior derivative estimates
\cite[Theorem~2.10]{GilbargTrudinger2001}
then yield
\[
\|w\|_{W^{2,p}(B_r)}
\le C\|w\|_{L^1(B_{r_1})}.
\]
Combining the estimates for $v$ and $w$
proves the assertion.
\end{proof}

\begin{lemma}\label{lem:Hartogs-order}
Let $E$ be a compact subset of a complex manifold. Suppose
$u_j,u$ are plurisubharmonic on a neighborhood of $E$, with
$u_j\to u$ locally in $L^1$ there. Then, for every $f\in C(E)$,
\[
\limsup_{j\to\infty}
\sup_{E}(u_j-f)
\le\sup_{E}(u-f).
\]
\end{lemma}

\begin{proof}
Suppose the assertion fails. Choose $A\in\R$ such that
\[
\sup_E(u-f)<A<
\limsup_{j\to\infty}\sup_E(u_j-f).
\]
After passing to a subsequence, there exist
$z_j\in E$ such that $z_j\to z\in E$ and
\[
u_j(z_j)-f(z_j)>A.
\]
Work in holomorphic coordinates near $z$. For all sufficiently
large $j$, take Euclidean balls of small fixed radius $r$ about
$z_j$ and $z$ contained in a relatively compact set of the chart. Since plurisubharmonic functions are subharmonic, the mean-value inequality \cite[Chapter~I, (4.12)(c) and (5.3)]{DemaillyCADG} gives
\[
u_j(z_j)\le\frac1{|B_r|}
\int_{B(z_j,r)}u_j\,dV.
\]
The $L^1$ convergence and $z_j\to z$ imply
\[
\int_{B(z_j,r)}u_j\,dV
\longrightarrow
\int_{B(z,r)}u\,dV.
\]
Indeed, this follows by writing $u_j=u+(u_j-u)$, using the $L^1$ convergence for the second term and the continuity of translations in $L^1$ for the first. Consequently,
\[
\limsup_{j\to\infty}u_j(z_j)
\le\frac1{|B_r|}\int_{B(z,r)}u\,dV.
\]
As $r\downarrow0$, the right-hand side converges to $u(z)$ by the mean-value inequality and upper semicontinuity of $u$. Thus
\[
\limsup_{j\to\infty}u_j(z_j)\le u(z).
\]
Since $f(z_j)\to f(z)$, we obtain
\[
A\le u(z)-f(z)\le\sup_E(u-f),
\]
contrary to the choice of $A$.
\end{proof}

\begin{theorem}\label{thm:rough-forward}
Let $T=\ddbar\psi\ge0$, where $\psi\in\PSH(\Omega)$, and suppose that
$\nu(T,z)=0$ for every $z\in\Omega$. Then there exists a smooth family
of K\"ahler forms $\omega_T(t)=t\theta_\mu+\ddbar U_T(t)$, $t>0$, with the following
properties:
\begin{enumerate}[label=\textup{(\roman*)}]
\item The potential $U_T$ satisfies
\[
\partial_tU_T=\log\frac{\omega_T(t)^n}{\mu},
\]
and
\begin{equation}\label{eq:rough-L1-trace}
U_T(t)\longrightarrow\psi
\quad\text{in }L^1_{\mathrm{loc}}(\Omega)
\quad(t\downarrow0).
\end{equation}
Consequently, $\omega_T(t)\rightharpoonup T$ in the sense of currents.

\item If $\widetilde U$ is a smooth solution of \eqref{eq:scalar-gauge}
on $\Omega\times(0,S]$, for some $S>0$,
such that $\widetilde U(t)\to\psi$ in $L^1_{\mathrm{loc}}(\Omega)$ as
$t\downarrow0$, then $\widetilde U\le U_T$ on $\Omega\times(0,S]$.

\item The metric flow depends only on $T$, independently of the smooth
approximations, the exhaustion, the global plurisubharmonic
potential, and the reference volume form.
\end{enumerate}
\end{theorem}

\begin{proof}
By \cite[Theorem~4.1]{HarveyLawsonPlis2016},
choose $\psi_j\in C^\infty(\Omega)\cap\SPSH(\Omega)$
with $\psi_j\downarrow\psi$. Monotone convergence
applied to $\psi_1-\psi_j$ gives
$\psi_j\to\psi$ in $L^1_{\mathrm{loc}}(\Omega)$.
Let $U_j$ and $\omega_j=t\theta+\ddbar U_j$ be the
solutions of Proposition~\ref{prop:fixed-domain-flow} on $\Omega_j$
with initial forms $\ddbar\psi_j$.

We first compare the solutions on successive exhaustion domains.
Since $\vartheta_j\ge0$, equation \eqref{eq:fixed-corner} gives
$\liminf_{t\downarrow0}\inf_{\Omega_j}(U_j-\psi_j)\ge0$.
Moreover, $U_{j+1}(t)\to\psi_{j+1}$ uniformly on
$\overline\Omega_j$, and $\psi_j\ge\psi_{j+1}$. Hence
\[
\liminf_{t\downarrow0}\inf_{\Omega_j}(U_j-U_{j+1})\ge0.
\]
By \eqref{eq:fixed-boundary-blow}, $U_j(t)$ tends to $+\infty$ at
$\partial\Omega_j$ locally uniformly in positive time, whereas
$U_{j+1}$ is bounded on $\overline\Omega_j$ over compact positive-time
intervals. Lemma~\ref{lem:scalar-comparison} therefore gives
\begin{equation}\label{eq:rough-monotone-m}
U_j\ge U_{j+1}\qquad\text{on }\Omega_j\times(0,\infty).
\end{equation}

The function $V_j(t)=\psi_j+Kt\rho+nt\log(Kt)-nt$
is a subsolution of \eqref{eq:scalar-gauge}, since
$t\theta+\ddbar V_j=\ddbar\psi_j+Kt\omega_E\ge Kt\omega_E>0$ and
\[
\partial_tV_j
=K\rho+n\log(Kt)
=\log\frac{(Kt\omega_E)^n}{\mu}
\le\log\frac{(t\theta+\ddbar V_j)^n}{\mu}.
\]
As $\psi_j$ and $\rho$ are bounded on $\overline\Omega_j$,
$V_j(t)\to\psi_j$ uniformly on $\Omega_j$ as $t\downarrow0$, and
$V_j$ is bounded near $\partial\Omega_j$ on compact positive-time
intervals. Lemma~\ref{lem:scalar-comparison} gives
\begin{equation}\label{eq:rough-lower}
U_j(t)\ge\psi_j+Kt\rho+nt\log(Kt)-nt.
\end{equation}
For $D\Subset\Omega$, choose $j_0$ with
$\overline D\Subset\Omega_{j_0}$. By
\eqref{eq:rough-monotone-m} and \eqref{eq:rough-lower},
\[
\psi+Kt\rho+nt\log(Kt)-nt\le U_j(t)\le U_{j_0}(t)
\]
on $D\times[\tau,S]$, for $j\ge j_0$ and $0<\tau<S<\infty$.
Since $\psi\in L^1_{\mathrm{loc}}(\Omega)$, we obtain
\[
\sup_{j\ge j_0}\sup_{\tau\le t\le S}
\|U_j(t)\|_{L^1(D)}<\infty.
\]

We next obtain volume-form estimates independent of $j$.
For every $\delta>0$, the time-shifted flow
$s\mapsto\omega_j(s+\delta)$ is complete on
$\Omega_j\times[0,\infty)$. By
\cite[Lemma~2.2]{L22},
\[
R_{\omega_j(t)}\ge-\frac{n}{t-\delta},
\qquad t>\delta.
\]
Letting $\delta\downarrow0$ gives
$R_{\omega_j(t)}\ge-n/t$.
Consequently,
\[
\partial_t^2(U_j-nt\log t)
=-R_{\omega_j(t)}-\frac nt\le0.
\]
Since $U_j(t)\to\psi_j$ locally uniformly as $t\downarrow0$, the
secant-slope inequality for the concave function
$t\mapsto U_j(t,z)-nt\log t$ gives
\begin{equation}\label{eq:rough-Ut-upper}
\partial_tU_j(t)\le\frac{U_j(t)-\psi_j}{t}+n.
\end{equation}
Fix $D\Subset D'\Subset\Omega$ and $0<\tau<S$, and choose $j_0$
with $\overline{D'}\Subset\Omega_{j_0}$. Adding a common constant to
$\psi$, $\psi_j$, and $U_j$, we may assume $\psi\le0$ on $D'$.
Since $\psi_j\ge\psi$ and $U_j\le U_{j_0}$ on $D'$ for $j\ge j_0$,
\eqref{eq:rough-Ut-upper} yields
\begin{equation}\label{eq:rough-density-upper}
\frac{\omega_j(t)^n}{\mu}
=e^{\partial_tU_j(t)}
\le C_{D',\tau,S}e^{-\psi/\tau}
\quad\text{on }D'\times[\tau,S].
\end{equation}

By Skoda's integrability theorem
\cite[(1.3) and \S2.1]{DiNezzaLu2017},
the vanishing of the Lelong numbers of
$T=\ddbar\psi$ implies
$e^{-A\psi}\in L^1_{\mathrm{loc}}(\Omega)$
for every $A>0$. Therefore, for
every $p>1$,
\begin{equation}\label{eq:rough-density-Lp}
\sup_{j\ge j_0}\sup_{\tau\le t\le S}
\left\|\frac{\omega_j(t)^n}{\mu}\right\|_{L^p(D')}
<\infty.
\end{equation}

Applying \cite[Lemma~2.1]{L22} to the same time-shifted complete
flows, with reference metric $h$, gives 
$\tr_{\omega_j(t)}h\le 2n/((n+1)\kappa(t-\delta))$ for
$t>\delta$. Letting $\delta\downarrow0$ yields
\begin{equation}\label{eq:rough-ambient-floor}
\omega_j(t)\ge a_nt\,h,\qquad a_n=\frac{(n+1)\kappa}{2n}.
\end{equation}

Fix concentric holomorphic coordinate balls
$D_0\Subset D_1\Subset D_2\Subset\Omega$ and $0<\tau<S$.
Choose $j_1$ with $\overline D_2\subset\Omega_{j_1}$. For
$j\ge j_1$, the eigenvalues of
$\omega_j(t)$ relative to $h$ on $D_2\times[\tau/2,2S]$ are at
least $a_n\tau/2$. Thus
\[
\tr_h\omega_j(t)\le n\left(\frac{2}{a_n\tau}\right)^{n-1}
\frac{\omega_j(t)^n}{h^n}.
\]
Since $\ddbar U_j=\omega_j-t\theta$, the real Euclidean
Laplacian of $U_j$ in this chart is bounded in $L^p(D_2)$ by the
trace estimate and
\eqref{eq:rough-density-Lp}, for every finite $p>1$.
The smooth term $t\theta$ has bounded coefficients on this
cylinder. The local $L^1$ bounds and Lemma~\ref{lem:W2p-L1} give
uniform $W^{2,p}(D_1)$ bounds for $U_j(t)$ on $[\tau/2,2S]$.
Sobolev embedding for $p>2n$ gives uniform $C^0$ bounds there.

For each $z\in D_1$, the function
$\sigma_{j,z}(t)=U_j(t,z)-nt\log t$ is concave.
For $\tau\le t\le S$, we have
\[
\frac{\sigma_{j,z}(2S)-\sigma_{j,z}(t)}{2S-t}
\le\sigma_{j,z}'(t)
\le\frac{\sigma_{j,z}(t)-\sigma_{j,z}(\tau/2)}
{t-\tau/2}.
\]
The uniform $C^0$ bounds on $[\tau/2,2S]$ and the inequalities
$t-\tau/2\ge\tau/2$, $2S-t\ge S$ give
\begin{equation}\label{eq:rough-Ut-two-sided}
|\partial_tU_j|\le C_{D_1,\tau,S}
\qquad\text{on }D_1\times[\tau,S].
\end{equation}
By \eqref{eq:scalar-gauge}, the volume ratio
$\omega_j(t)^n/\mu$ is now bounded above and below
by positive constants. Combining these bounds with
\eqref{eq:rough-ambient-floor}, we obtain
\begin{equation}\label{eq:rough-two-sided}
C^{-1}h\le\omega_j(t)\le Ch
\qquad\text{on }D_1\times[\tau,S],
\end{equation}
where $C$ is independent of $j$.

Applying \eqref{eq:rough-two-sided} with $\tau$ replaced by $\tau/2$ gives uniform metric equivalence on $D_1\times[\tau/2,S]$. Hence \cite[Corollary~1.2]{ShermanWeinkove2012}, applied on coordinate balls to the time-shifted flows $s\mapsto\omega_j(s+\tau/2)$, gives uniform interior derivative bounds on $D_0\times[\tau,S]$. The flow equation
\eqref{eq:KRF} gives the corresponding time-derivative bounds.
Taking the flat trace of $\ddbar U_j=\omega_j-t\theta$ and
using the local bounds for $U_j$, the interior elliptic estimates
\cite[\S6.3.1, Theorem~2]{PDEbookEvans_1998}, applied in
Euclidean coordinates, give bounds for all spatial derivatives
of $U_j$. Differentiating \eqref{eq:scalar-gauge} gives the corresponding mixed space--time derivative bounds.

As in the proof of Proposition~\ref{prop:fixed-domain-flow},
the uniform interior estimates give
$C^\infty_{\mathrm{loc}}$ precompactness.
By \eqref{eq:rough-monotone-m}, the pointwise limit
is unique, so the entire sequence converges:
\[
U_j\longrightarrow U_T
\quad\text{in }C^\infty_{\mathrm{loc}}(\Omega\times(0,\infty)).
\]
Passing to the limit in \eqref{eq:scalar-gauge} and
\eqref{eq:rough-two-sided} shows that
$\omega_T(t)=t\theta+\ddbar U_T(t)>0$ is a smooth K\"ahler--Ricci flow.

We now verify \eqref{eq:rough-L1-trace}. Passing to the limit in \eqref{eq:rough-lower} gives
\begin{equation}\label{eq:rough-limit-lower}
U_T(t)\ge\psi+Kt\rho+nt\log(Kt)-nt.
\end{equation}
On every compact $E\subset\Omega$, \eqref{eq:rough-limit-lower}
gives $(\psi-U_T(t))_+\to0$ in $L^1(E)$. For any $j_0$ with
$E\subset\Omega_{j_0}$, monotonicity gives $U_T\le U_{j_0}$ on $E$.
The locally uniform trace of $U_{j_0}$ therefore gives
\[
(U_T(t)-\psi)_+\le\psi_{j_0}-\psi+o_{j_0}(1)
\quad\text{on }E
\]
almost everywhere, with the error uniform on $E$ as $t\downarrow0$.
Letting $t\downarrow0$ and then $j_0\to\infty$ proves the local
$L^1$ trace. Applying $\ddbar$ distributionally and using
$t\theta\to0$ gives $\omega_T(t)\rightharpoonup T$.

It remains to prove \textup{(ii)} and \textup{(iii)}. Let $\widetilde U$ satisfy the hypotheses of \textup{(ii)}.
Fix $j$, $t_0\in(0,S]$, and $\gamma>0$. On each coordinate chart $\widetilde U(t)+tb_z$ is psh and has the
same local $L^1$ trace $\psi$. Applying Lemma~\ref{lem:Hartogs-order}
on finitely many compact chart sets covering $\overline\Omega_j$,
using $\psi\le\psi_j$ and the uniform smallness of $tb_z$, gives
\[
\widetilde U(\delta)\le\psi_j+\frac{\gamma}{2}
\qquad\text{on }\overline\Omega_j
\]
for all sufficiently small $\delta>0$. By
\eqref{eq:fixed-corner}, after decreasing $\delta$ if necessary,
\[
U_j(\delta)\ge\psi_j-\frac{\gamma}{2}
\qquad\text{on }\Omega_j.
\]
Thus $U_j(\delta)+\gamma\ge\widetilde U(\delta)$.
For $t\in[\delta,t_0]$, equation \eqref{eq:fixed-boundary-blow}
gives $U_j(t)\to+\infty$ at $\partial\Omega_j$, locally uniformly
in time, whereas $\widetilde U$ is bounded on
$\overline\Omega_j\times[\delta,t_0]$. Applying
Lemma~\ref{lem:scalar-comparison} on $[\delta,t_0]$ gives $U_j(t_0)+\gamma\ge\widetilde U(t_0)$.
Letting $\gamma\downarrow0$ and then $j\to\infty$, we obtain
\begin{equation}\label{eq:rough-maximal}
U_T\ge\widetilde U.
\end{equation}

We now prove \textup{(iii)}. If $U_T$ and $\widetilde U_T$
are obtained from different approximations or exhaustions
with the same initial potential $\psi$, applying
\eqref{eq:rough-maximal} to each gives
$U_T\ge\widetilde U_T$ and
$\widetilde U_T\ge U_T$. Thus $U_T=\widetilde U_T$.

If $\psi'\in\PSH(\Omega)$ is another global potential
for $T$, then $\ddbar(\psi'-\psi)=0$ distributionally.
By Weyl's lemma
\cite[Theorem~IX.25]{ReedSimonII1975},
there exists a smooth real pluriharmonic function
$p$ such that $\psi'=\psi+p$.
Let $U_T'$ be the solution constructed with
initial potential $\psi'$. Since $p$ is
pluriharmonic and independent of time,
$U_T+p$ and $U_T'-p$ satisfy
\eqref{eq:scalar-gauge} with initial traces
$\psi'$ and $\psi$, respectively. Applying
\eqref{eq:rough-maximal} gives
\[
U_T'\ge U_T+p,
\qquad U_T\ge U_T'-p.
\]
Hence $U_T'=U_T+p$, and the metric is unchanged.
For $\widehat\mu=e^b\mu$, the transformation
$U\mapsto U-tb$ preserves the order and gives a bijection of scalar
solutions with trace $\psi$ in the two gauges. It therefore takes
the maximal solution to the maximal solution. This also proves
existence for any smooth positive reference volume and independence
of the metric from that choice.
\end{proof}

We next establish the transformation law of
the canonical flow under biholomorphisms and
constant rescaling. Its proof uses the maximality
property \eqref{eq:rough-maximal}.

\begin{proposition}\label{prop:functoriality}
Let $\Omega\Subset Z$ and $\Omega'\Subset Z'$
be Stein domains in complex manifolds, and suppose
that neighborhoods of their closures admit
holomorphic immersions into Euclidean spaces.
Let $F:\Omega'\to\Omega$ be a biholomorphism.
Suppose $T=\ddbar\psi\ge0$ on $\Omega$,
where $\psi\in\PSH(\Omega)$ and
$\nu(T,z)=0$ for every $z\in\Omega$.
Then, for every $c>0$,
\begin{equation}\label{eq:functoriality}
\omega_{cF^*T}(t)=cF^*\omega_T(t/c),\qquad t>0.
\end{equation}
\end{proposition}

\begin{proof}
Let $\mu,\mu'$ be the reference volumes on $\Omega,\Omega'$,
respectively. The ratio $J_F=F^*\mu/\mu'$ satisfies
$\ddbar\log J_F=F^*\theta_\mu-\theta_{\mu'}$. 
Define
\[
\widetilde U(t)
=cU_T(t/c)\circ F
+t(n\log c+\log J_F).
\]
Since
$\ddbar\log J_F=F^*\theta_\mu-\theta_{\mu'}$,
we have
\[
t\theta_{\mu'}+\ddbar\widetilde U(t)
=cF^*\omega_T(t/c).
\]
Moreover,
\[
\frac{(cF^*\omega_T(t/c))^n}{\mu'}
=c^nJ_F
\left(\frac{\omega_T(t/c)^n}{\mu}\right)\circ F.
\]
Consequently, by \eqref{eq:scalar-gauge},
\[
\begin{aligned}
\partial_t\widetilde U
&=(\partial_sU_T)(t/c)\circ F
+n\log c+\log J_F =\log\frac{
(t\theta_{\mu'}+\ddbar\widetilde U(t))^n
}{\mu'}.
\end{aligned}
\]
Its initial trace in $L^1_{\mathrm{loc}}(\Omega')$ is $c\psi\circ F$.
By biholomorphic invariance and scaling of Lelong numbers,
Theorem~\ref{thm:rough-forward} applies to $cF^*T$. 
By \eqref{eq:rough-maximal},
$\widetilde U\le U_{cF^*T}$.
Conversely, define on $\Omega$
\[
\widehat U(s)
=c^{-1}U_{cF^*T}(cs)\circ F^{-1}
+s\bigl(-n\log c+\log J_{F^{-1}}\bigr).
\]
The same calculation shows that $\widehat U$
satisfies \eqref{eq:scalar-gauge} with reference
volume $\mu$ and initial trace $\psi$.

Applying \eqref{eq:rough-maximal} again,
we obtain $\widehat U\le U_T$.
Substituting $s=t/c$ and composing with $F$,
we obtain
\[
\begin{aligned}
U_{cF^*T}(t)
&\le cU_T(t/c)\circ F
+t(n\log c+\log J_F)\\
&=\widetilde U(t),
\end{aligned}
\]
where we used $J_{F^{-1}}\circ F=J_F^{-1}$.
Thus $U_{cF^*T}=\widetilde U$, and taking
the associated K\"ahler forms proves
\eqref{eq:functoriality}.
\end{proof}

\subsection{Instantaneous completeness}
\label{sec:completeness}

Although each $\omega_j(t)$ is complete on $\Omega_j$,
the smooth local convergence $\omega_j\to\omega_T$
does not imply completeness of $\omega_T(t)$.
We therefore use the flows $\omega_{j,\varepsilon}(t)$ of
Proposition~\ref{prop:fixed-domain-flow}, with initial metric
\eqref{eq:fixed-regularized-initial} for $\eta=\ddbar\psi_j$.
We prove \eqref{eq:regularized-volume-floor} and
\eqref{eq:regularized-loggrad} with bounds independent of
$j$ and $\varepsilon$.

\begin{lemma}\label{lem:regularized-volume-floor}
For every $j$, $0<\varepsilon<1$, and $t>0$,
\begin{equation}\label{eq:regularized-volume-floor}
\omega_{j,\varepsilon}(t)^n\ge (Kt)^n\omega_E^n
\qquad\text{on }\Omega_j.
\end{equation}
\end{lemma}

\begin{proof}
By \cite[Lemma~2.2]{L22}, the complete flows
$\omega_{j,\varepsilon}(t)$ of
Proposition~\ref{prop:fixed-domain-flow} satisfy
\begin{equation}\label{eq:regularized-scalar-lower}
R_{\omega_{j,\varepsilon}(t)}\ge-\frac nt.
\end{equation}
Fix $t>0$, write $\omega=\omega_{j,\varepsilon}(t)$ and set
$r=\omega_E^n/\omega^n$. By \eqref{eq:baromega-det} and
\eqref{eq:fixed-regularized-equivalence}, $\omega^n/\mu$ diverges
at $\partial\Omega_j$. Since $\omega_E$ is smooth on a neighborhood
of $\overline\Omega_j$ in $\Omega$, $r$ tends to zero there.
At an interior maximum point of $r$, the identity
$\Ric_{\omega_E}=-K\omega_E$ from
\eqref{eq:KE-normalization} and the
arithmetic--geometric mean inequality give
\[
0\ge\Delta_\omega\log r
=R_\omega+K\tr_\omega\omega_E
\ge -\frac nt+nKr^{1/n}.
\]
Since $r\le(Kt)^{-n}$ at its maximum point, we obtain
\eqref{eq:regularized-volume-floor}.
\end{proof}

Let $g_0$ be a complete Riemannian metric on $X$.
By interpolating between $g_\beta$ near $\overline\Omega$
and $g_0$ outside a compact subset of $X$, choose a
Riemannian metric $g_*$ such that $g_*=g_\beta$ near
$\overline\Omega$ and $g_*\ge c_0g_0$ for some $c_0>0$.
In particular, $g_*$ is complete. Let
$d(z)=\dist_{g_*}(z,\partial_X\Omega)$. We prove the boundary
volume estimate
\begin{equation}\label{eq:MokYau-density-blow}
\frac{\omega_E^n}{\mu}(z)
\ge \frac{c}{d(z)^2(-\log d(z))^2}
\end{equation}
for sufficiently small $d(z)>0$. We adapt the boundary-volume comparison of
Mok--Yau~\cite[\S3.2]{MokYau}
to the inclusion $\Omega\hookrightarrow(X,\beta)$.

By Hopf--Rinow, every $z\in\Omega$ can be joined to
$\partial_X\Omega$ by a minimizing $g_*$-geodesic. Its open segment
lies in $\Omega$, since an earlier exit would be a closer boundary
point. On this segment the metric is $g_\beta$. Conversely, any
finite-length $\beta$-geodesic in $\Omega$ leaving every compact
subset is a Cauchy curve for $g_*$ and converges to $\partial_X\Omega$.
Thus $d$ is the exit distance for the inclusion
$\Omega\hookrightarrow(X,\beta)$.
For a Stein coordinate ball $W\subset X$, let
$\Delta_X=\{(z,z):z\in X\}$ be the diagonal.
Since $X$ is Hausdorff, the graph
$\Gamma_W=\Delta_X\cap(\Omega\times W)$ is a closed complex
submanifold of the Stein manifold $\Omega\times W$.
Since $\Gamma_W$ is biholomorphic to $\Omega\cap W$,
the latter is Stein. Thus $\Omega$ is locally
pseudoconvex over $X$, and
\cite[Theorem~2.1.11]{Dingoyan1999OkaLevi} gives
$\ddbar(-\log d)\ge-C\beta$ as currents on $\{d<d_*\}$ for
some $C,d_*>0$. Choose $0<d_0<\min(d_*,e^{-2})$, and set
\[
u=\max\{-\log d,-\log d_0\},\qquad L(v)=2v-2\log v\quad(v>2),
\qquad \nu=e^{A\varphi+L(u)}\mu,\quad A\ge2C+1.
\]
On the open set $\{d<d_*\}$, both
$-\log d+C\varphi$ and $-\log d_0+C\varphi$
are plurisubharmonic. Their maximum is
$u+C\varphi$ there. On the open set
$\{d>d_0\}$, this function equals
$-\log d_0+C\varphi$ and is plurisubharmonic.
Since $d_0<d_*$, these open sets cover
$\Omega$. Thus $u+C\varphi$ is continuous
and plurisubharmonic on $\Omega$.

By Lemma~\ref{lem:stein-KE-exhaustion} and \eqref{eq:baromega-det},
$\nu/\eta_j^n$ tends to zero as $z\to\partial\Omega_j$,
so it attains a positive maximum at some $x\in\Omega_j$. Put
$M=\log(\nu/\eta_j^n)(x)$. By maximality,
\[
L(u)\le\log(\eta_j^n/\mu)+M-A\varphi
\]
on $\Omega_j$, with equality at $x$. Since $L'>0$ on $(2,\infty)$,
the function
\[
\widetilde u=L^{-1}\bigl(\log(\eta_j^n/\mu)+M-A\varphi\bigr)
\]
is smooth near $x$ and satisfies $\widetilde u\ge u$,
with equality at $x$. Since $u+C\varphi$ is
plurisubharmonic, the mean-value inequality for
subharmonic functions on complex lines gives
$\ddbar\widetilde u(x)\ge-C\beta(x)$.
Differentiating the definition of $\widetilde u$ above and
using $\Ric_{\eta_j}=-K\eta_j$ from
Lemma~\ref{lem:stein-KE-exhaustion}, $L'<2$, and
$L''(u)=2/u^2$, we obtain at $x$
\[
K\eta_j
=\theta+A\beta+L'(u)\ddbar\widetilde u+
        \frac{2}{u^2}\ii\,\partial\widetilde u\wedge\bar\partial\widetilde u
\ge\beta+\frac{2}{u^2}\ii\,\partial\widetilde u\wedge\bar\partial\widetilde u.
\]
If $d(x)<d_0$, then $e^{-\widetilde u}$ touches $d$ from below.
Choose $y$ a short positive distance from $x$ along a minimizing
$g_*$-geodesic from $x$ to $\partial_X\Omega$. The triangle inequality gives
$d(y)=d(x)-\dist_{g_*}(x,y)$, and $z\mapsto d(y)+\dist_{g_*}(z,y)$
touches $d$ from above at $x$.
For $y$ sufficiently close to $x$, the function
$z\mapsto\dist_{g_*}(z,y)$ is smooth near $x$ and has
real gradient of norm one. Since $e^{-\widetilde u}$ and
$z\mapsto d(y)+\dist_{g_*}(z,y)$ touch $d$ at $x$
from below and above, respectively, their differentials agree there.
Consequently,
$|\partial\widetilde u|_\beta^2=1/(2d(x)^2)$.
Here the real metric associated to $\beta$ is $2\operatorname{Re}(\beta_{i\bar j}dz^i d\bar z^j)$,
so that $|\partial f|_\beta^2=\tfrac12|df|_{g_\beta}^2$.
The rank-one determinant formula yields at $x$
\[
(K\eta_j)^n\ge \left(1+\frac1{u^2d(x)^2}\right)\beta^n
\ge e^{L(u)}\beta^n.
\]
If $d(x)\ge d_0$, then $u(x)=-\log d_0$ and
$K\eta_j\ge\beta$ implies
\[
(K\eta_j)^n(x)\ge\beta^n(x)
=e^{-L(-\log d_0)}e^{L(u(x))}\beta^n(x).
\]
Since $e^{A\varphi}\mu/\beta^n$
is bounded on $\overline\Omega$, in either case
$\eta_j^n(x)\ge c_1\nu(x)$ with $c_1>0$ independent of $j$.
Since $\nu/\eta_j^n$ attains its maximum at $x$,
we have $\eta_j^n\ge c_1\nu$ on $\Omega_j$.
Letting $j\to\infty$ and using Lemma~\ref{lem:stein-KE-exhaustion}
proves \eqref{eq:MokYau-density-blow}. When $Z=\C^n$, we may take $\iota$ to be
the identity, so that $\beta=\omega_{\mathrm{Euc}}$.
The distance $d$ then agrees with Euclidean
boundary distance up to the fixed normalization
of the associated Riemannian metric. Since
$h^n/dV_{\mathrm E}$ is bounded above and
below on $\overline\Omega$, the estimate
\eqref{eq:MokYau-density-blow} holds with
Euclidean boundary distance and volume.

If $\omega(t)$ satisfies \eqref{eq:KRF}, then in holomorphic
coordinates the function $F_i=\log(\omega(t)^n/dV_i)$ satisfies
$\partial_tF_i=\Delta_{\omega(t)}F_i=-R_{\omega(t)}$. Thus
\begin{equation}\label{eq:F-caloric}
(\partial_t-\Delta_\omega)F_i=0.
\end{equation}

\begin{proposition}\label{prop:regularized-loggrad}
Let $\mu$ be a smooth positive volume on a neighborhood of
$\overline\Omega$. Set $F_E=\log(\omega_E^n/\mu)$,
$m_E=\inf_\Omega F_E$, and
$F_{j,\varepsilon}(t)=\log(\omega_{j,\varepsilon}(t)^n/\mu)$.
Then $m_E>-\infty$. For every $\tau>0$, set
$A_\tau=n\log(K\tau/2)+m_E-1$. There is $C_\tau<\infty$ independent
of $j,\varepsilon$ such that
\begin{equation}\label{eq:regularized-loggrad}
\left|\partial\log\bigl(F_{j,\varepsilon}(\tau)-A_\tau\bigr)
\right|_{\omega_{j,\varepsilon}(\tau)}^2\le C_\tau.
\end{equation}
When $\Omega\subset\C^n$ and $\mu=dV_{\mathrm E}$, the estimate holds
with $C_\tau=8n/\tau$.
\end{proposition}

\begin{proof}
Since $\mu/h^n$ is smooth and positive on a
neighborhood of $\overline\Omega$,
\eqref{eq:MokYau-density-blow} implies that
$F_E$ tends to $+\infty$ at
$\partial_X\Omega$ and is bounded below. Fix $j,\varepsilon,\tau$
and write $\omega=\omega_{j,\varepsilon}$ and $F=F_{j,\varepsilon}$.
By \eqref{eq:regularized-volume-floor}, $v=F-A_\tau\ge1$ on
$\Omega_j\times I$, where $I=[\tau/2,\tau]$.
By \cite[Lemma~2.1]{L22},
\begin{equation}\label{eq:regularized-metric-floor}
\omega(t)\ge a_n t h,\qquad t\in I.
\end{equation}

Choose finitely many holomorphic coordinate charts covering
$\overline\Omega$, and relatively compact smaller sets $V_i$ in
them still covering $\overline\Omega$. Take $0\le\zeta_i\le1$
smooth, equal to one on $V_i$ and compactly supported in
the corresponding chart. Write $\mu=e^{b_i}dV_i$ there.
Choose $B\ge0$ such that $|b_i|\le B$ on a neighborhood
of $\operatorname{supp}\zeta_i$ in the $i$-th coordinate chart,
for every $i$. On these neighborhoods set $v_i=v+b_i+B$,
$w_i=\log v_i$, and $q_i=|\partial w_i|_\omega^2$.
Then $1\le v\le v_i\le v+2B$ there. For $P=\partial_t-\Delta_\omega$,
\eqref{eq:F-caloric} gives $Pv_i=0$ and
\begin{equation}\label{eq:w-P}
Pw_i=q_i,\qquad (w_i)_t=-R_\omega/v_i\le 2n/\tau .
\end{equation}
The inequality $(w_i)_t\le2n/\tau$ follows from
\eqref{eq:regularized-scalar-lower} and $v_i\ge1$.
For a smooth real-valued function $w$, set
$q_w=|\partial w|_\omega^2=\tfrac12|\nabla^g w|_g^2$.
Writing the Bochner formula
\cite[\S2, (2.2)]{GangLiu2014} in real notation
and using $\partial_tg=-\Ric_g$ from \eqref{eq:KRF},
\[
\begin{aligned}
\partial_tq_w
&=2\operatorname{Re}
\langle\partial w_t,\partial w\rangle_\omega
+\tfrac12\Ric_g(\nabla^g w,\nabla^g w),\\
\Delta_\omega q_w
&=2\operatorname{Re}
\langle\partial\Delta_\omega w,\partial w\rangle_\omega
+\tfrac12\Ric_g(\nabla^g w,\nabla^g w)
+\tfrac12|\nabla_g^2w|_g^2.
\end{aligned}
\]
Subtracting gives
\begin{equation}\label{eq:K\"ahler-bochner}
P|\partial w|_\omega^2
=2\operatorname{Re}\langle\partial(Pw),\partial w\rangle_\omega
-\tfrac12|\nabla_g^2w|_g^2.
\end{equation} 
By Cauchy--Schwarz and $\Delta_g=2\Delta_\omega$,
\[
\tfrac12|\nabla_g^2w|_g^2
\ge\frac{1}{4n}(\Delta_g w)^2
=\frac1n(\Delta_\omega w)^2.
\]
Consequently,
\begin{equation}\label{eq:q-P}
Pq_i\le
2\operatorname{Re}\langle\partial q_i,\partial w_i\rangle_\omega
-\frac1n(\Delta_\omega w_i)^2 .
\end{equation}

Fix $j$ and $\tau>0$. By the bounded-geometry
construction above, every point of $\Omega_j$
is the center of a holomorphic coordinate ball
whose radius and reference metric bounds are
uniform for this fixed $j$.
By \eqref{eq:fixed-regularized-equivalence},
the metrics $\omega=\omega_{j,\varepsilon}$
are uniformly equivalent to the flat metric
on each of these balls throughout
$[\tau/4,\tau]$, independently of
$\varepsilon$ and of the chart center.
Applying
\cite[Theorem~1.1(i)]{ShermanWeinkove2012}
on a smaller concentric ball, after shifting
time by $\tau/4$, and using the uniform
first derivatives of $\overline\omega_j$
in these coordinates, gives
\begin{equation}\label{eq:fixed-positive-time-derivative}
|\nabla_{\overline\omega_j}\omega|_{\overline\omega_j}
\le C_{j,\tau}
\qquad\text{on }\Omega_j\times I.
\end{equation}
The constant is independent of
$\varepsilon$ and of the chart center. Differentiating the determinant
and using \eqref{eq:fixed-regularized-equivalence} gives
\begin{equation}\label{eq:fixed-logdet-derivative}
\left|\partial\log(\omega^n/\overline\omega_j^n)\right|_\omega
\le C_{j,\tau}\qquad\text{on }\Omega_j\times I.
\end{equation}

By \eqref{eq:baromega-det}, in the $i$-th
coordinate chart we have
\[
F_i
=\log\frac{\omega^n}{\overline\omega_j^n}
+a\vartheta_j+\ell_j
+\log\frac{h^n}{dV_i}.
\]
Since $\ell_j\in C^\infty(\overline\Omega_j)$
and $\log(h^n/dV_i)$ is smooth in the
$i$-th coordinate chart, these functions
and their first derivatives are bounded on
$\operatorname{supp}\zeta_i\cap\overline\Omega_j$.
Using \eqref{eq:fixed-regularized-equivalence},
\eqref{eq:fixed-logdet-derivative}, and
$|\partial\vartheta_j|_{\overline\omega_j}\le1$,
we obtain
\[
F_i=a\vartheta_j+O_{j,\tau}(1),
\qquad
|\partial F_i|_\omega\le C_{j,\tau}
\]
on $\operatorname{supp}\zeta_i\cap\Omega_j$,
uniformly in $\varepsilon$ and $t\in I$.
Since $\vartheta_j(z)\to+\infty$ as $z\to\partial\Omega_j$
by \eqref{eq:background-potential}, it follows that
\begin{equation}\label{eq:fixed-boundary-q-decay}
q_i=|\partial F_i|_\omega^2/v_i^2
\longrightarrow0
\qquad(z\to\partial\Omega_j)
\end{equation}
uniformly for $t\in I$ and
$z\in\operatorname{supp}\zeta_i\cap\Omega_j$.

Put $\xi_i=\zeta_i^2$. By
\eqref{eq:regularized-metric-floor},
there exists $L_\tau>0$, independent of
$j,\varepsilon,i$, such that
\[
|\Delta_\omega\xi_i|\le L_\tau,\qquad
|\partial\xi_i|_\omega^2/\xi_i=4|\partial\zeta_i|_\omega^2
\le L_\tau\quad(\xi_i>0).
\]
Write $s=t-\tau/2$ and $Q_i=s\xi_i q_i$.
Extend $Q_i$ by zero outside $\operatorname{supp}\zeta_i$.
By \eqref{eq:fixed-boundary-q-decay}, this extension tends to zero
as $z\to\partial\Omega_j$, uniformly for $t\in I$;
moreover, $Q_i(\tau/2,\cdot)=0$.
If $Q_i\equiv0$, the estimate on $V_i$ is immediate.
Otherwise, $Q_i$ attains a positive maximum at a point
$(z_i,t_i)\in\Omega_j\times(\tau/2,\tau]$ with $\xi_i(z_i)>0$.
At $(z_i,t_i)$, we have $PQ_i\ge0$, using the left time
derivative if $t_i=\tau$, and
$\partial q_i=-q_i\partial\xi_i/\xi_i$.
Writing $Q,\xi,q,w$ for $Q_i,\xi_i,q_i,w_i$,
the product rule and \eqref{eq:q-P} give
\begin{equation}\label{eq:loggrad-max-inequality}
0\le PQ\le
\xi q+3s L_\tau q+2s q\sqrt{L_\tau\xi q}
-\frac{s\xi}{n}(\Delta_\omega w)^2 .
\end{equation}
At $(z_i,t_i)$, the terms
$-sq\Delta_\omega\xi-2s\operatorname{Re}
\langle\partial\xi,\partial q\rangle_\omega$ are at most
$3sL_\tau q$, while
$2s\xi\operatorname{Re}\langle\partial q,\partial w\rangle_\omega
\le 2sq\sqrt{L_\tau\xi q}$.
If $q\le4n/\tau$ at $(z_i,t_i)$, then $Q\le2n$. Otherwise
$\Delta_\omega w=w_t-q<-q/2$ by \eqref{eq:w-P}. Dividing
\eqref{eq:loggrad-max-inequality} by $q$ and applying Young's
inequality yields
\[
\frac{Q}{4n}\le1+3sL_\tau+2\sqrt{sL_\tau Q}
\le1+(3+8n)sL_\tau+\frac{Q}{8n}.
\]
Thus $Q\le8n[1+(3+8n)(\tau/2)L_\tau]$.
Since $\xi_i=1$ on $V_i$ and $s=\tau/2$ at $t=\tau$,
\[
q_i(\tau)\le\frac{16n}{\tau}
\left[1+(3+8n)\frac{\tau L_\tau}{2}\right]
\qquad\text{on }V_i\cap\Omega_j.
\]
Moreover,
\[
\partial\log v=(v_i/v)\partial w_i-\partial b_i/v.
\]
By $1\le v\le v_i\le v+2B$ and
\eqref{eq:regularized-metric-floor}, the bounds for
$|\partial b_i|_h$ on $V_i$ yield
\eqref{eq:regularized-loggrad}, since the sets $V_i$
cover $\overline\Omega$.

When $\mu=dV_{\mathrm E}$, the function $F$
satisfies \eqref{eq:F-caloric} globally.
We may therefore apply the preceding argument
with $v_i=v$, $\xi_i=1$, and $L_\tau=0$
on $\Omega_j$.
By \eqref{eq:fixed-logdet-derivative} and
\eqref{eq:baromega-det},
$|\partial F|_\omega\le C_{j,\tau}$ and
$F=a\vartheta_j+O_{j,\tau}(1)$ on $\Omega_j\times I$,
since
\[
\overline\omega_j^n
=e^{a\vartheta_j+\widetilde\ell_j}dV_{\mathrm E},
\qquad
\widetilde\ell_j\in C^\infty(\overline\Omega_j).
\]
Hence $q=|\partial F|_\omega^2/v^2\to0$ as
$z\to\partial\Omega_j$, uniformly for $t\in I$.

If $Q\equiv0$, the estimate is immediate.
Otherwise, $Q$ attains a positive maximum
at $(z_0,t_0)\in\Omega_j\times(\tau/2,\tau]$.
At this point, if
$q\le4n/\tau$, then $Q\le2n$.
Otherwise, by \eqref{eq:w-P},
$\Delta_\omega w=w_t-q<-q/2$.
Thus \eqref{eq:loggrad-max-inequality},
with $L_\tau=0$, gives
\[
0\le q-\frac{sq^2}{4n},
\]
and hence $Q\le4n$. In either case,
$Q\le4n$, so at $t=\tau$,
\[
|\partial\log v|_\omega^2
\le\frac{8n}{\tau}.
\]
\end{proof}

\begin{theorem}\label{thm:instant-complete}
Let $T=\ddbar\psi\ge0$, where $\psi\in\PSH(\Omega)$ and
$\nu(T,z)=0$ for every $z\in\Omega$. The canonical flow
$\omega_T(t)$ from Theorem~\ref{thm:rough-forward} is complete
for every $t>0$.
\end{theorem}

\begin{proof}
Fix $\tau>0$. Pass to $\varepsilon_k\downarrow0$ using
\eqref{eq:fixed-flow-smooth-limit}, and then let $j\to\infty$
using the $C^\infty_{\mathrm{loc}}$ convergence in the proof of
Theorem~\ref{thm:rough-forward}. Equations
\eqref{eq:regularized-volume-floor} and
\eqref{eq:regularized-loggrad} give
\begin{align}
\omega_T(\tau)^n&\ge(K\tau)^n\omega_E^n,\label{eq:limit-volume-floor}\\
|\partial\log(F_T(\tau)-A_\tau)|_{\omega_T(\tau)}^2
&\le C_\tau,\label{eq:limit-loggrad}
\end{align}
where $F_T(t)=\log(\omega_T(t)^n/\mu)$, and
$C_\tau=8n/\tau$ when $\mu=dV_{\mathrm E}$.
By \eqref{eq:regularized-volume-floor},
\[
F_{j,\varepsilon}(\tau)-A_\tau
\ge 1+n\log2+F_E-m_E\ge1+n\log2.
\]
This inequality also holds for $\omega_j(\tau)$ and
$\omega_T(\tau)$ by smooth local convergence, so the
logarithms in \eqref{eq:limit-loggrad} are defined.

Set $\Xi_\tau=\log(F_T(\tau)-A_\tau)$. By
\eqref{eq:limit-volume-floor} and \eqref{eq:MokYau-density-blow},
$\Xi_\tau(z)\to+\infty$ as $z\to\partial_X\Omega$.
Hence $\{z\in\Omega:\Xi_\tau(z)\le c\}$ is compact in $\Omega$
for every $c\in\R$. Let $g_\tau$ be the Riemannian metric
associated to $\omega_T(\tau)$. By \eqref{eq:limit-loggrad}
and the convention in Section~\ref{sec:preliminaries},
$|\nabla^{g_\tau}\Xi_\tau|_{g_\tau}^2
=2|\partial\Xi_\tau|_{\omega_T(\tau)}^2\le2C_\tau$. Integrating along
piecewise smooth curves and taking the infimum of their lengths gives
\[
|\Xi_\tau(x)-\Xi_\tau(y)|
\le \sqrt{2C_\tau}\,d_{g_\tau}(x,y).
\]
When $\mu=dV_{\mathrm E}$, the constant is $4\sqrt{n/\tau}$.
If $(z_k)$ is Cauchy with respect to $d_{g_\tau}$,
then $\sup_kd_{g_\tau}(z_k,z_1)<\infty$.
The preceding inequality gives $\sup_k\Xi_\tau(z_k)<\infty$,
so $(z_k)$ lies in a compact sublevel set of $\Xi_\tau$.
It has a subsequence converging in the Riemannian topology.
The Cauchy property gives convergence of the sequence, and
Hopf--Rinow gives completeness.
\end{proof}

\subsection{Self-similarity and the soliton equation}

We now prove Theorem~\ref{thm:current-to-soliton}.
The transformation law \eqref{eq:functoriality},
together with \eqref{eq:current-homogeneity},
gives the following self-similarity identity.

\begin{proposition}\label{prop:homogeneous-selfsimilar}
Let $T\in\mathcal T_{Y,K}(\Omega)$.
Then the K\"ahler--Ricci flow $\omega_T(t)$
constructed in Theorem~\ref{thm:rough-forward}
satisfies
\begin{equation}\label{eq:canonical-selfsimilar}
\omega_T(t)=e^{-Ks}\Phi_s^*\omega_T(e^{Ks}t),
\qquad s\in\R,\quad t>0.
\end{equation}
In particular, $\omega=\omega_T(K^{-1})$ is
a complete K\"ahler metric satisfying
\eqref{eq:expanding-soliton}, and $\omega_T(t)$
is the associated self-similar flow
\eqref{eq:selfsimilar-flow}.
\end{proposition}

\begin{proof}
By Propositions~\ref{prop:current-exact-psh}
and \ref{prop:current-zero-lelong}, the current
$T$ admits a global plurisubharmonic potential
and has vanishing Lelong numbers. Thus
Theorem~\ref{thm:rough-forward} provides
the flow $\omega_T(t)$, which is complete
for every $t>0$ by Theorem~\ref{thm:instant-complete}.

For every $s\in\R$,
\eqref{eq:current-homogeneity} gives
$e^{-Ks}\Phi_s^*T=T$.
Applying \eqref{eq:functoriality} with
$F=\Phi_s$ and $c=e^{-Ks}$ proves
\eqref{eq:canonical-selfsimilar}.
Differentiating this identity with respect
to $s$ at $s=0$ and using \eqref{eq:KRF},
we obtain
\[
0=-K\omega_T(t)+\Lie_Y\omega_T(t)
-Kt\,\Ric_{\omega_T(t)}.
\]
At $t=K^{-1}$, this is precisely
\eqref{eq:expanding-soliton} for
$\omega=\omega_T(K^{-1})$.
Finally, substituting
$s=-K^{-1}\log(Kt)$ into
\eqref{eq:canonical-selfsimilar}
recovers \eqref{eq:selfsimilar-flow}.
\end{proof}

\begin{proof}[Proof of Theorem~\ref{thm:current-to-soliton}]
Let $T\ne0$ be the current in the statement
of the theorem. Proposition~\ref{prop:homogeneous-selfsimilar}
gives a complete expanding soliton
$\omega=\omega_T(K^{-1})$ with soliton
vector field $Y$ and expansion constant $K$.

Suppose that $\Lie_Y\omega=0$.
Since $\frac{d}{ds}\Phi_s^*\omega
=\Phi_s^*(\Lie_Y\omega)=0$
and $\Phi_0=\id_\Omega$, we have
$\Phi_s^*\omega=\omega$ for every $s\in\R$,
using the completeness of $Y$.
Equation~\eqref{eq:selfsimilar-flow} then gives
$\omega_T(t)=Kt\,\omega$, which converges
weakly to zero as $t\downarrow0$.
On the other hand, \eqref{eq:rough-L1-trace}
gives $\omega_T(t)\rightharpoonup T\ne0$,
a contradiction. Thus
$\Lie_Y\omega\not\equiv0$, proving
that the soliton is nontrivial.
\end{proof}

\subsection{Products with Stein domains}
\label{sec:product-expander}
We conclude with the proof of
Theorem~\ref{thm:product-stability}.

\begin{proof}[Proof of Theorem~\ref{thm:product-stability}]
Suppose that $D\times D'$ admits a complete expanding
K\"ahler--Ricci soliton with complete soliton vector field.
The product of the given immersions is a holomorphic
immersion on a neighborhood of $\overline{D\times D'}$.
By Remarks~\ref{rem:stein-consequences}\textup{(2)},
the product is Stein. Let $\chi$ be a smooth
strictly plurisubharmonic exhaustion of $D\times D'$,
and fix $z_0\in D$. Since $\{z_0\}\times D'$ is
closed in $D\times D'$, the restriction
$\chi(z_0,\cdot)$ is a smooth strictly
plurisubharmonic exhaustion of $D'$.
Thus $D'$ is Stein.

Conversely, suppose $D'$ is Stein. If $D'$ is a point, the
conclusion follows from the soliton on $D$ and constant rescaling.
Otherwise let $(\omega_D,Y_D)$ be the complete nontrivial expanding
soliton on $D$ with expansion constant $K_0$.
By Remarks~\ref{rem:stein-consequences}\textup{(2)},
$D'$ admits a complete K\"ahler--Einstein
metric $\omega_{D'}$. After a constant rescaling,
we may assume that
\[
\Ric_{\omega_{D'}}=-K_0\omega_{D'}.
\]
Let $\operatorname{pr}_D$ and $\operatorname{pr}_{D'}$
denote the projections from $D\times D'$ onto its factors,
and set
\[
\omega=\operatorname{pr}_D^*\omega_D
       +\operatorname{pr}_{D'}^*\omega_{D'},
\qquad Y=(Y_D,0).
\]
The metric $\omega$ is complete as the Riemannian product
of two complete metrics. The vector field $Y$ is complete
and real holomorphic, with flow
$(z,w)\mapsto(\Phi_s^D(z),w)$, where $\Phi_s^D$ is the
flow of $Y_D$. In product coordinates, the coefficient matrix
of $\omega$ is block diagonal, so its determinant is the
product of the determinants of the two factors. Consequently,
using \eqref{eq:expanding-soliton} on $D$, we obtain
\[
\begin{aligned}
\Ric_\omega+K_0\omega
&=\operatorname{pr}_D^*(\Ric_{\omega_D}+K_0\omega_D)
 +\operatorname{pr}_{D'}^*(\Ric_{\omega_{D'}}+K_0\omega_{D'})\\
&=\operatorname{pr}_D^*(\Lie_{Y_D}\omega_D)
 =\Lie_Y\omega\not\equiv0.
\end{aligned}
\]
Hence $(\omega,Y)$ is a complete nontrivial expanding
soliton on $D\times D'$ with expansion constant $K_0$.

For any $K>0$, define
$\widehat\omega=(K_0/K)\omega$ and
$\widehat Y=(K/K_0)Y=((K/K_0)Y_D,0)$.
Since constant rescaling leaves the Ricci form unchanged
and $\Lie_{\widehat Y}\widehat\omega=\Lie_Y\omega$,
the pair $(\widehat\omega,\widehat Y)$ satisfies
\eqref{eq:expanding-soliton} with expansion constant $K$.
Both the metric and the vector field remain complete,
and $\Lie_{\widehat Y}\widehat\omega\not\equiv0$.
\end{proof}

Since $D$ admits a complete expanding soliton,
Remarks~\ref{rem:stein-consequences}\textup{(2)}
imply that $D$ is Stein. When $D'$ is
also Stein, those remarks give complete
K\"ahler--Einstein metrics on both
factors with the same negative Einstein constant.
Their product therefore defines a complete
trivial expanding soliton on $D\times D'$.


\section{Uniformization and Expanding Solitons}
\label{sec:boundary-uniformization}

To apply the uniformization theorems, we use
the following consequence of Kaup's properness results
\cite[S\"atze~1.2--1.3, pp.~306--308]{Kaup1968}.
The statement is formulated for complex spaces, as it
will also be used in Section~\ref{sec:bordered-stein-spaces}.

\begin{lemma}\label{lem:bounded-coordinate-action}
Let $D$ be a reduced irreducible complex space whose bounded
holomorphic functions separate points. Then $\Aut(D)$ acts properly.
If $D$ is an open subset of a complex space $X$ with
$D\Subset X$, every subgroup of $\Aut(D)$ with noncompact closure has
an orbit accumulating at $\partial_X D$.
\end{lemma}

\begin{proof}
Since bounded holomorphic functions separate points
of $D$, every subset of $D$ containing at least two
points admits a bounded holomorphic function whose
restriction to that subset is nonconstant. Thus $D$ is
boundedly $K$-complete in the sense of Kaup
\cite[p.~306]{Kaup1968}. As observed there, bounded
$K$-completeness is preserved under coverings. In particular,
the universal covering of $D$ is boundedly $K$-complete.
Kaup's Satz~1.2~\cite[p.~306]{Kaup1968} then shows that $D$
is weakly hyperbolic. Since $D$ is irreducible, Satz~1.3
\cite[p.~308]{Kaup1968} implies that $\Aut(D)$ is a real Lie
group acting properly on $D$.

Fix $p\in D$. By properness of the action, the map
\[
\Aut(D)\times D\longrightarrow D\times D,\qquad
(F,z)\longmapsto (z,F(z))
\]
is proper. Consequently, for every compact set $K\subset D$, the set
\[
\{F\in\Aut(D):F(p)\in K\}
\]
is compact, since its product with $\{p\}$ is the inverse image
of $\{p\}\times K$ under this map. Hence the orbit map
$F\mapsto F(p)$ is proper.

Now suppose that $D\Subset X$, and let $H\subset\Aut(D)$
be a subgroup with noncompact closure. Choose an increasing
compact exhaustion $\{K_j\}$ of $D$ with
$K_j\subset\operatorname{int}K_{j+1}$. Since the orbit map
is proper, each set
\[
\{F\in\Aut(D):F(p)\in K_j\}
\]
is compact and therefore cannot contain $H$. We may thus
choose $F_j\in H$ with $F_j(p)\notin K_j$. 
Since every compact subset of $D$ is contained in some $K_J$,
the points $F_j(p)\notin K_j$ eventually lie outside any
given compact subset of $D$. On the other hand,
$\overline D$ is compact, so a subsequence of
$\{F_j(p)\}$ converges to a point $q\in\overline D$.
Necessarily $q\notin D$, and hence $q\in\partial_XD$.
\end{proof}

The Siegel half-space
\[
\mathcal S=\{(z',w)\in\C^{n-1}\times\C:
2\Re w+|z'|^2<0\}
\]
and, in dimension two, the domains $\mathcal G_H$ of
\eqref{eq:surface-model} are invariant under translations
in the imaginary $w$-direction. More generally, let
$G\subset\C^{n-1}\times\C$ be a domain invariant under
the translations
\[
\tau_s(z',w)=(z',w+\ii s),\qquad s\in\R.
\]
Writing $w=u+\ii v$, these automorphisms form a free
$\R$-action generated by the complete real holomorphic
vector field $Y_0=\partial/\partial v$. The action is also
proper. Indeed, if $(p_j,\tau_{s_j}(p_j))$ remains in a
compact subset of $G\times G$, then both $p_j$ and
$\tau_{s_j}(p_j)$ remain in compact subsets of $G$, while
\[
s_j=\Im w(\tau_{s_j}(p_j))-\Im w(p_j)
\]
is bounded. After passing to a subsequence, both $p_j$
and $s_j$ converge, which proves properness.

For every $K>0$, the forms
\begin{equation}\label{eq:translation-currents}
T_\pm=\ddbar e^{\pm K\Im w}
=\frac{K^2}{4}e^{\pm K\Im w}\,\ii\,dw\wedge d\bar w
\end{equation}
are nonzero smooth semipositive closed $(1,1)$-forms on $G$.
They satisfy
\[
\tau_s^*T_+=e^{Ks}T_+,
\qquad
\tau_{-s}^*T_-=e^{Ks}T_-,
\qquad s\in\R.
\]
Thus $T_+$ and $T_-$ satisfy
\eqref{eq:current-homogeneity} for the complete vector
fields $Y_0$ and $-Y_0$, respectively.


\subsection{Strongly pseudoconvex domains}

\begin{proof}[Proof of Theorem~\ref{thm:smooth-ball}]
Suppose that $\Omega$ admits a complete nontrivial
expanding soliton with complete real holomorphic soliton
field. By Lemma~\ref{pcvx:lem:compact-closure-rigidity},
its one-parameter automorphism group has noncompact
closure in $\Aut(\Omega)$. Since $\Omega$ is a connected
complex manifold, it is reduced and irreducible. Since bounded holomorphic functions on $\Omega$
separate points, Lemma~\ref{lem:bounded-coordinate-action}
implies there exist $p\in\Omega$, $F_j\in\Aut(\Omega)$,
and $\zeta\in\partial_X\Omega$ such that
$F_j(p)\to\zeta$.

We assumed $\zeta$ is a $C^2$ strongly
pseudoconvex boundary point. In suitable local
holomorphic coordinates the boundary is strongly
convex near $\zeta$, so Theorem~1 of
Gaussier--Kim--Krantz
\cite[p.~762]{GaussierKimKrantz2002}
gives $\Omega\cong\B^n$.

Conversely, suppose that $\Omega$ is biholomorphic
to $\B^n$. The Cayley map
\[
C(\zeta',\zeta_n)=
\left(\frac{\sqrt{2}\zeta'}{1-\zeta_n},
      -\frac{1+\zeta_n}{1-\zeta_n}\right)
\]
identifies $\B^n$ biholomorphically with the
Siegel half-space $\mathcal S$ introduced above. For every $K>0$, the currents $C^*T_\pm$ from
\eqref{eq:translation-currents} are nonzero,
positive and closed, and satisfy
\eqref{eq:current-homogeneity} for the
complete real holomorphic vector fields
$\pm C^*Y_0$, respectively.

Since $\B^n\Subset\C^n$ is Stein,
Theorem~\ref{thm:current-to-soliton}
gives complete nontrivial expanding solitons
on $\B^n$ with expansion constant $K$
and these prescribed fields. Pulling them
back by any biholomorphism
$F:\Omega\to\B^n$ gives the corresponding
solitons on $\Omega$, with fields
$\pm(C\circ F)^*Y_0$.
\end{proof}

\subsection{Finite-type domains in complex surfaces}

\begin{corollary}\label{thm:surface-local}
Let $\Omega\Subset X$ be a domain in a complex surface $X$
such that a neighborhood of $\overline\Omega$ admits a
holomorphic embedding into $\C^N$.
Let $\zeta\in\partial_X\Omega$ be a boundary orbit
accumulation point. Assume that the boundary is
$C^\infty$ smooth, pseudoconvex and of finite type
near $\zeta$, with type $2m$ at $\zeta$.
Then $\Omega$ is biholomorphic to a domain
$\mathcal G_H$ of the form \eqref{eq:surface-model},
where $H$ is a nonzero real homogeneous subharmonic
polynomial of degree $2m$ without harmonic terms.
For every biholomorphism $F:\Omega\to\mathcal G_H$
and every $K>0$, there exist complete nontrivial
expanding K\"ahler--Ricci solitons on $\Omega$
with expansion constant $K$ and soliton vector
fields $F^*Y_0$ and $-F^*Y_0$.
\end{corollary}

\begin{proof}
By Seshadri--Verma
\cite[Theorem~1.1]{SeshadriVerma2008},
$\Omega$ is biholomorphic to a domain $\mathcal G_H$
with $H$ as stated. We first observe that
$\mathcal G_H$ is Stein. Indeed, writing
$r=2\Re w+H(z,\bar z)$, the function
\[
|z|^2+|w|^2+\max\{0,-\log(-r)\}
\]
is a continuous plurisubharmonic exhaustion
of $\mathcal G_H$.

Fix a biholomorphism $F:\Omega\to\mathcal G_H$
and $K>0$. The currents $F^*T_+$ and $F^*T_-$,
defined in \eqref{eq:translation-currents},
are nonzero positive closed $(1,1)$-currents
satisfying \eqref{eq:current-homogeneity}
with respect to $F^*Y_0$ and $-F^*Y_0$,
respectively. Since $\Omega$ is Stein,
Theorem~\ref{thm:current-to-soliton}
applies to these currents and gives the
required complete nontrivial expanding solitons.
\end{proof}

\begin{corollary}\label{thm:surface-global}
Let $\Omega\Subset X$ be a domain in a complex surface.
Suppose that a neighborhood of $\overline\Omega$
admits a holomorphic embedding into $\C^N$.
Assume that every boundary orbit accumulation point
of $\Aut(\Omega)$ has a neighborhood in which
the boundary is $C^\infty$ smooth, pseudoconvex,
and of finite type.
Then the following are equivalent:
\begin{enumerate}[label=\textup{(\roman*)}]
\item $\Aut(\Omega)$ is noncompact.
\item $\Aut_0(\Omega)$ is noncompact.
\item $\Omega$ is biholomorphic to a model \eqref{eq:surface-model}
with $H$ as in Corollary~\ref{thm:surface-local}.
\item $\Omega$ admits a complete nontrivial expanding
K\"ahler--Ricci soliton with complete real holomorphic soliton field.
\end{enumerate}
\end{corollary}

\begin{proof}
We prove
\textup{(i)}$\Rightarrow$\textup{(iii)}
$\Rightarrow$\textup{(iv)}
$\Rightarrow$\textup{(ii)}
$\Rightarrow$\textup{(i)}.

Suppose \textup{(i)} holds. The coordinate functions of
the holomorphic embedding of a neighborhood of
$\overline\Omega$ are bounded on $\Omega$ and separate
points. By Lemma~\ref{lem:bounded-coordinate-action},
there exists a boundary orbit accumulation point
$\zeta\in\partial_X\Omega$. The boundary is
$C^\infty$ smooth, pseudoconvex and of finite type
near $\zeta$ by hypothesis. Hence
Corollary~\ref{thm:surface-local} gives \textup{(iii)}.

Suppose next that \textup{(iii)} holds, and let
$F:\Omega\to\mathcal G_H$ be a biholomorphism.
Since $\mathcal G_H$ is Stein, so is $\Omega$.
For any $K>0$, the nonzero positive closed current
$F^*T_+$ of \eqref{eq:translation-currents}
satisfies \eqref{eq:current-homogeneity} for the
complete real holomorphic vector field $F^*Y_0$.
Theorem~\ref{thm:current-to-soliton}
therefore gives \textup{(iv)}.

Now suppose \textup{(iv)} holds. If
$\Aut_0(\Omega)$ were compact, the flow of the
complete soliton vector field would have compact
closure in $\Aut(\Omega)$. By
Lemma~\ref{pcvx:lem:compact-closure-rigidity},
the soliton would be trivial, a contradiction.
Thus \textup{(ii)} holds.

Finally, \textup{(ii)} implies \textup{(i)},
since $\Aut_0(\Omega)$ is a subgroup of
$\Aut(\Omega)$.
\end{proof}

\section{Automorphism Groups of Regular Loci of Stein Domains}
\label{sec:bordered-stein-spaces}

Let $Z$ be a reduced Stein space of pure dimension
$n\ge2$, and let $D\Subset Z$ be a connected Stein
domain. Since $D$ is open in $Z$, we have
\[
M:=\Reg D=D\cap\Reg Z,\qquad
S:=\Sing D=D\cap\Sing Z.
\]
We first establish the normalization and embedding
properties used to prove Theorem~\ref{thm:compact}.
We then derive the soliton obstruction and prove
the two uniformization corollaries.

\begin{lemma}\label{lem:normalization}
Suppose that $S$ is finite. Let
$\nu_Z:\widehat Z\to Z$ be the normalization
of $Z$, and set $\widehat D=\nu_Z^{-1}(D)$.
Then $\widehat Z$ and $\widehat D$ are Stein,
$\nu_Z:\widehat D\to D$ is the normalization
of $D$, and $\widehat D$ is a finite disjoint
union of irreducible normal Stein spaces.
Moreover,
\begin{equation}\label{eq:punctured-normalization}
M\simeq\widehat D\setminus E,
\qquad E=\nu_Z^{-1}(S),
\end{equation}
where $E$ is finite. If $S$ is nonempty, every
connected component of $\widehat D$ meets $E$.
There exists a relatively compact Stein
neighborhood $W$ of $\overline{\widehat D}$
in $\widehat Z$ admitting a proper holomorphic
embedding into $\C^N$ for some $N$.
\end{lemma}

\begin{proof}
By the decomposition lemma for analytic sets
\cite[Chapter~A, \S3.5]
{GrauertRemmert2004}, the irreducible components
$\{D_\alpha\}_{\alpha\in A}$ of $D$ form a
locally finite family of analytic subsets of $D$,
each of which is closed in $D$.
At every point of $\Reg D$, the local ring is
isomorphic to an algebra of convergent power
series, and hence the corresponding germ is
irreducible
\cite[Chapter~A, \S3.7--3.8, p.~21]
{GrauertRemmert2004}.
Since $D$ is pure $n$-dimensional, two
irreducible components meeting at a regular
point would have the same local branch there.
By the identity theorem for analytic sets
\cite[Chapter~A, \S3.5, p.~19]
{GrauertRemmert2004}, they would then coincide.
Thus distinct irreducible components of $D$
can intersect only in $S$.
Suppose that $S$ is nonempty. If some
$D_\alpha$ were disjoint from $S$, it would
also be disjoint from every other component.
By local finiteness, $D_\alpha$ would then
be both open and closed in $D$, contradicting
connectedness. Thus every component meets $S$. Since the family is locally finite, each point of $S$ belongs
to only finitely many components. As $S$ is
finite, $D$ has only finitely many irreducible
components. If $S$ is empty, $D$ is a connected
complex manifold and hence irreducible.

Normalization is finite, commutes with restriction
to open subsets, and is biholomorphic over the
regular locus
\cite[Definition~1.1.94 and Theorem~1.1.95]
{GreuelLossenShustin2025}.
It follows that $\widehat D$ is a finite
disjoint union of irreducible normal components,
each mapping onto an irreducible component of
$D$. Since every component of $D$ meets $S$
when $S$ is nonempty, each component of
$\widehat D$ then meets $E$. The normalizations
of the Stein spaces $Z$ and $D$ are Stein by
\cite[Theorem~1]{Narasimhan1962}.
Since $S$ is finite and $\nu_Z$ is finite,
$E$ is finite, and restriction gives
\eqref{eq:punctured-normalization}.
Since the normalization map $\nu_Z$ is finite,
it is proper
\cite[Chapter~A, \S3.8, and
Chapter~I, \S1.1]{GrauertRemmert2004}.
As $\overline D$ is compact and
\[
\overline{\widehat D}
\subseteq \nu_Z^{-1}(\overline D),
\]
it follows that $\overline{\widehat D}$
is compact in $\widehat Z$. Set $L=\overline{\widehat D}$, and let
$\widehat L$ be its holomorphic hull in
$\widehat Z$. Since $\widehat Z$ is Stein,
$\widehat L$ is compact. Choose an open set
$V$ with $\widehat L\subset V\Subset\widehat Z$.
For every $p\in\partial V$, the definition
of the holomorphic hull gives
$f_p\in\mathcal O(\widehat Z)$ such that,
after multiplication by a positive constant,
\[
\sup_L|f_p|<1<|f_p(p)|.
\]
By compactness of $\partial V$, finitely many
such functions $f_1,\ldots,f_k$ suffice to define
\[
U=\bigcap_{\ell=1}^k
\{x\in\widehat Z:|f_\ell(x)|<1\}
\]
with $L\subset U$ and
$\overline U\cap\partial V=\varnothing$.

The graph of $(f_1,\ldots,f_k)$ over $U$
is a closed complex subspace of
$\widehat Z\times\unitdisc^k$. Since this
product is Stein, so is $U$
\cite[Chapter~V, \S1, Theorem~1]
{GrauertRemmert2004}.
Let $W$ be the union of those connected components
of $U$ which meet $L$. By compactness of $L$,
only finitely many components occur. Each is
contained in $V$, since it meets $V$ and its
closure avoids $\partial V$. Thus $W$ is a
relatively compact Stein neighborhood of $L$.
Since $W\Subset\widehat Z$, a finite covering
of $\overline W$ by local embedding charts
gives a uniform bound on the local embedding
dimension of $W$. By Cartan's theorem~A
\cite[Chapter~IV, Fundamental Theorem]{GrauertRemmert2004},
the local embedding coordinates can be chosen
as germs of global holomorphic functions on
the Stein space $W$. Thus $W$ is locally of
bounded type in the sense of Narasimhan, and
\cite[Theorem~6, p.~931]{NarasimhanEmbedding}
gives a proper holomorphic embedding
$W\hookrightarrow\C^N$ for some $N$.
\end{proof}

\begin{proof}[Proof of Theorem~\ref{thm:compact}]
By Lemma~\ref{lem:normalization},
$\widehat D$ is a finite disjoint union of
irreducible normal Stein spaces
$D_1,\ldots,D_r$, with
$E_i=E\cap D_i$ finite and nonempty, and
$M\simeq\widehat D\setminus E$.
Let $j:W\hookrightarrow\C^N$ be the proper
holomorphic embedding from that lemma.
The coordinate functions of $j$ are bounded
on $\widehat D$, since
$\overline{\widehat D}$ is compact in $W$.

Let $f\in\Aut(M)$. As $\widehat D$ is normal
and $E$ has codimension $n\ge2$, the coordinate
functions of $j\circ f$ extend across $E$
by the Riemann extension theorem
\cite[Theorem~1.1.98]{GreuelLossenShustin2025}.
They define a holomorphic map
$H:\widehat D\to\C^N$. Since $M$ is dense
in $\widehat D$ and $H(M)=j(M)$, continuity gives
\[
H(\widehat D)\subset j(\overline{\widehat D})
\subset j(W).
\]
Since $j(W)$ is a closed complex subspace of
$\C^N$ and $\widehat D$ is reduced, $H$
factors holomorphically through $j(W)$.
Thus $F=j^{-1}\circ H$ is a holomorphic map
from $\widehat D$ to $W$ with
\[
F|_M=f,\qquad F(\widehat D)
\subset\overline{\widehat D}.
\]

We claim that $F(\widehat D)\subset\widehat D$.
Otherwise, there is $a\in E$ such that
$b=F(a)\in\partial_W\widehat D$.
Since $F(M)=M$, the fiber $F^{-1}(b)$ is
contained in $E$, so $a$ is isolated in
this fiber. Since $E$ is finite, apply \cite[Theorem~1.1.66]{GreuelLossenShustin2025}
to the restriction of $F$ to a neighborhood
of $a$ disjoint from $E\setminus\{a\}$.
This gives neighborhoods $U_a$ of $a$ and
$V_b$ of $b$ such that $F:U_a\to V_b$ is
finite and $U_a\cap E=\{a\}$.
Since $W$ is normal, we may shrink $V_b$
to a connected irreducible neighborhood,
replacing $U_a$ by $U_a\cap F^{-1}(V_b)$. \cite[Corollary~1.1.68]{GreuelLossenShustin2025}
shows that $F(U_a)$ is a closed analytic subset
of $V_b$. Since $f$ is a biholomorphism of
$M$, the set $f(U_a\setminus\{a\})$ is
nonempty and open in $V_b$. As $V_b$ is
irreducible, it follows that $F(U_a)=V_b$.
Since $F$ maps $U_a\setminus\{a\}$ into $M$,
\[
V_b\setminus\{b\}\subset\widehat D.
\]

Choose a relatively compact connected open
neighborhood $P$ of $b$ in $V_b$ with compact
frontier $K=\partial_{V_b}P$.
Then $K\subset\widehat D$. For every
$g\in\mathcal O(\widehat D)$, the restriction
of $g$ to $V_b\setminus\{b\}$ extends across
$b$ by the Riemann extension theorem
\cite[Theorem~1.1.98]{GreuelLossenShustin2025}.
The maximum-modulus principle on the
irreducible complex space $P$
\cite[Introduction, p.~XIX]{GrauertRemmert2004}
gives
\[
|g(x)|\le\sup_K|g|,
\qquad x\in P\setminus\{b\}.
\]
Consequently, the holomorphic hull of $K$
in $\widehat D$ contains $P\setminus\{b\}$,
which is not relatively compact in
$\widehat D$ because it accumulates at
$b\notin\widehat D$. This contradicts
holomorphic convexity of the Stein space
$\widehat D$. Hence
$F(\widehat D)\subset\widehat D$.

Applying the same argument to $f^{-1}$
gives a holomorphic self-map of $\widehat D$
extending $f^{-1}$. The compositions equal the
identity on the dense set $M$, and therefore
on $\widehat D$. Thus $F$ is an automorphism
preserving $E$, and restriction is a continuous
surjective homomorphism
\begin{equation}\label{eq:group-iso}
\Aut(\widehat D;E)\longrightarrow\Aut(M),
\qquad
\Aut(\widehat D;E)
=\{G\in\Aut(\widehat D):G(E)=E\}.
\end{equation}
Continuity follows directly from the compact-open
topology: a compact subset of $M$ is compact in
$\widehat D$, and an open subset of $M$ is open
in $\widehat D$.

Finally, the coordinate functions of $j$
are bounded on each $D_i$ and separate points.
Since $D_i$ is irreducible,
Lemma~\ref{lem:bounded-coordinate-action}
shows that $\Aut(D_i)$ acts properly on $D_i$.
Fix $a_i\in E_i$. The evaluation map
$\operatorname{ev}_{a_i}:\Aut(D_i)\to D_i$
is proper, and hence
$\operatorname{ev}_{a_i}^{-1}(E_i)$ is compact.
Moreover,
\[
\Aut(D_i;E_i)
=\bigcap_{e\in E_i}\operatorname{ev}_e^{-1}(E_i)
\]
is a closed subset of that compact set.
Thus $\Aut(D_i;E_i)$ is compact.
The subgroup of $\Aut(\widehat D;E)$ preserving
all components $D_i$ is the finite product
$\prod_{i=1}^r\Aut(D_i;E_i)$. It has finite
index, since automorphisms permute the finitely
many components. Consequently,
$\Aut(\widehat D;E)$ is compact, and so is
its continuous image $\Aut(M)$.
\end{proof}

\begin{proof}[Proof of Corollary~\ref{cor:obstruction}]
Suppose that $D$ is singular and that $M$ admits
such a soliton $(\omega,Y)$. By
Theorem~\ref{thm:compact}, $\Aut(M)$ is compact.
The flow of $Y$ preserves each connected component
of $M$, and its restriction has compact closure
in the automorphism group of that component.
Lemma~\ref{pcvx:lem:compact-closure-rigidity}
therefore gives $\Lie_Y\omega=0$ on every
component of $M$, contradicting nontriviality.
\end{proof}

Theorem~\ref{thm:compact} also rules out nonzero
homogeneous positive closed $(1,1)$-currents on $M$
when $S$ is finite and nonempty. Indeed, if $Y$ is a complete
real holomorphic vector field on $M$, compactness
of $\Aut(M)$ gives a sequence $s_j\to+\infty$
with $\Phi_{s_j}\to\id_M$, as in the proof of
Lemma~\ref{pcvx:lem:compact-closure-rigidity}.
Consequently, $\Phi_{s_j}^*T\rightharpoonup T$
for every current $T$ on $M$. If
$T\in\mathcal T_{Y,K}(M)$ with $K>0$ and $T\ne0$,
choose a compactly supported test form $\beta$
with $T(\beta)\ne0$. Homogeneity gives
\[
(\Phi_{s_j}^*T)(\beta)=e^{Ks_j}T(\beta),
\]
contradicting the weak convergence. Thus $T=0$.

We finish with the ball and finite-type surface uniformization corollaries.

\begin{proof}[Proof of Corollary~\ref{cor:ball}]
Suppose that $M$ admits a complete nontrivial
expanding soliton. By Corollary~\ref{cor:obstruction},
$D$ is smooth. Let $\pi:\widetilde Z\to Z$
be a resolution of singularities that is
biholomorphic over $\Reg Z$
\cite{BierstoneMilman1997}, and set
$\widetilde D=\pi^{-1}(D)$.
Then $\pi:\widetilde D\to D$ is biholomorphic.
Since $\pi$ is proper and $D\Subset Z$,
we have $\widetilde D\Subset\widetilde Z$.

As $Z$ is Stein, the restrictions to $D$ of
its holomorphic functions are bounded and
separate points. The same therefore holds
on $\widetilde D$. Every boundary orbit
accumulation point of $\Aut(\widetilde D)$
projects under $\pi$ to a boundary orbit
accumulation point of $\Aut(D)$. By hypothesis,
the latter belongs to $\Reg Z$ and is a $C^2$
strongly pseudoconvex boundary point.
Since $\pi$ is biholomorphic near such points,
the corresponding boundary points of
$\widetilde D$ have the same property.
Theorem~\ref{thm:smooth-ball}, applied to
$\widetilde D\Subset\widetilde Z$, now gives
$D\cong\B^n$.

Conversely, suppose that $F:D\to\B^n$ is a
biholomorphism. The proof of
Theorem~\ref{thm:smooth-ball} establishes the
existence of a complete nontrivial expanding
soliton $(\omega_{\B},Y_{\B})$ on $\B^n$.
Then $(F^*\omega_{\B},F^*Y_{\B})$ is a
complete nontrivial expanding soliton on $D=M$.
\end{proof}

\begin{proof}[Proof of Corollary~\ref{thm:bordered-surface}]
Since $\partial_ZD\subset\Reg Z$,
$\Sing D=D\cap\Sing Z$ is a compact analytic
subset of the Stein space $D$, and hence finite
\cite[p.~XIX]{GrauertRemmert2004}.
If $D$ is smooth, then $\overline D\subset\Reg Z$.
In this case, let
$W$ and $j:W\hookrightarrow\C^N$ be as in
Lemma~\ref{lem:normalization}. Since $\nu_Z$ is
biholomorphic over $\Reg Z$, the embedding $j$
induces a holomorphic embedding into $\C^N$
of the neighborhood
$V=\nu_Z(W\cap\nu_Z^{-1}(\Reg Z))$ of
$\overline D$ in $\Reg Z$.

Suppose that $M$ admits a complete nontrivial
expanding soliton. Corollary~\ref{cor:obstruction}
implies that $D$ is smooth. The boundary
hypothesis and Corollary~\ref{thm:surface-global},
applied with ambient surface $\Reg Z$, give
a biholomorphism $D\cong\mathcal G_H$.

Conversely, suppose that $F:D\to\mathcal G_H$
is a biholomorphism. Then $D$ is smooth, so
the preceding embedding is available. For any
$K>0$, the currents
$F^*T_+$ and $F^*T_-$ of
\eqref{eq:translation-currents} are nonzero,
positive and closed, and satisfy
\eqref{eq:current-homogeneity} for the complete
real holomorphic vector fields $F^*Y_0$
and $-F^*Y_0$, respectively. Since $D$ is
Stein, Theorem~\ref{thm:current-to-soliton}
applied to these currents gives complete
nontrivial expanding solitons with the
prescribed fields and expansion constant $K$.
\end{proof}

\section*{Acknowledgments}
The author thanks Max Hallgren and Yuan Yuan
for helpful discussions concerning the subject
of this paper. The author is also grateful to
Peter Ebenfelt and Ming Xiao for their
encouragement and support.
\section*{AI-use declaration}
The author used OpenAI's ChatGPT to assist with
locating relevant literature, translating passages
from non-English mathematical sources, and improving
the grammar, organization, and exposition of the
manuscript. The author takes full responsibility
for the contents of this paper.

\bibliographystyle{amsplain}
\bibliography{refKRsoliton}

\end{document}